\documentclass{amsart}

\usepackage{amsmath,amssymb,amsbsy,amsthm,amsfonts,mathtools,dsfont,hyperref,mathrsfs}

\usepackage{enumitem}
\usepackage{wrapfig}
\usepackage[top=3.5cm,bottom=3.5cm,left=3.5cm,right=3.5cm]{geometry}
\usepackage[normalem]{ulem}

\usepackage{graphicx}
\usepackage{tikz}
\usepackage{tikzscale}
\usetikzlibrary{intersections,decorations.markings}
\tikzset{>=latex}
\usepackage{pgfplots}
\pgfplotsset{compat=newest}
\usepgfplotslibrary{fillbetween}
\usepackage[normalem]{ulem}

\allowdisplaybreaks[3]

\newcommand{\N}{\mathbb{N}}
\newcommand{\R}{\mathbb{R}}

\newcommand{\tr}{\operatorname{Tr}}

\begin{document}

\newcommand{\ddd}{\,{\rm d}}
%\numberwithin{equation}{section} \marginparwidth=2cm

\def\note#1{\marginpar{\small #1}}
\def\tens#1{\pmb{\mathsf{#1}}}
\def\vec#1{\boldsymbol{#1}}
\def\norm#1{\left|\!\left| #1 \right|\!\right|}
\def\fnorm#1{|\!| #1 |\!|}
\def\abs#1{\left| #1 \right|}
\def\ti{\text{I}}
\def\tii{\text{I\!I}}
\def\tiii{\text{I\!I\!I}}

\newcommand{\loc}{{\rm loc}}
\def\diver{\mathop{\mathrm{div}}\nolimits}
\def\grad{\mathop{\mathrm{grad}}\nolimits}
\def\Div{\mathop{\mathrm{Div}}\nolimits}
\def\Grad{\mathop{\mathrm{Grad}}\nolimits}
\def\cof{\mathop{\mathrm{cof}}\nolimits}
\def\det{\mathop{\mathrm{det}}\nolimits}
\def\lin{\mathop{\mathrm{span}}\nolimits}
\def\pr{\noindent \textbf{Proof: }}

\def\pp#1#2{\frac{\partial #1}{\partial #2}}
\def\dd#1#2{\frac{\d #1}{\d #2}}
\def\vec#1{\boldsymbol{#1}}

\def\0{\vec{0}}
\def\A{\mathcal{A}}
\def\B{\mathcal{B}}
\def\b{\vec{b}}
\def\C{\mathcal{C}}
\def\c{\vec{c}}
\def\D{\vec{Dv}}
\def\DD{\vec{D}}
\def\BB{\vec{B}}
\def\e{\varepsilon}
\def\er{\epsilon}
\def\f{\vec{f}}
\def\F{\vec{F}}
\def\tF{\tilde{\F}}
\def\g{\vec{g}}
\def\G{\vec{G}}
\def\cG{\mathcal{\G}}
\def\H{\vec{H}}
\def\cH{\mathcal{H}}
\def\I{\vec{I}}
\def\Im{\text{Im}}
\def\j{\vec{j}}
\def\J{\vec{J}}
\def\dd{\vec{d}}
\def\k{\vec{k}}
\def\n{\vec{n}}
\def\q{\vec{q}}
\def\S{\vec{S}}
\def\s{\vec{s}}
\def\T{\vec{T}}
\def\u{\vec{u}}
\def\vp{\vec{\varphi}}
\def\vv{\vec{v}}
\def\vvt{\vv_{\tau}}
\def\vov{\vv\otimes\vv}
\def\cV{\mathcal{V}}
\def\w{\vec{w}}
\def\W{\vec{W}}
\def\x{\vec{x}}
\def\z{\vec{Z}}
\def\tz{\tilde{\z}}
\def\Z{\vec{Z}}
\def\X{\vec{X}}
\def\Y{\vec{Y}}
\def\balfa{\vec{\alpha}}

\def\Ge{\G_{\e}}
\def\ge{\g_{\e}}
\def\fidv{\phi_{\delta}(|\vv|^2)}
\def\fidve{\phi_{\delta}(|\ve|^2)}
\def\fidvd{\phi_{\delta}(|\vd|^2)}

\def\Ae{\A_\e}
\def\Aee{\Ae^\e}
\def\Aeetilde{\tilde{A}_\e^\e}
\def\Be{\B_\e}
\def\Bee{\Be^\e}
\def\De{\DD\ve}
\def\DDe{\DD^\e}
\def\Dvd{\DD\vd}
\def\oD{\overline{\DD}}
\def\tD{\tilde{\DD}}
\def\Dn{\DD^\e}
\def\Dno{\overline{\Dn}}
\def\Dnt{\tilde{\Dn}}
\def\Dm{\DD^\eta}
\def\Dmo{\overline{\Dm}}
\def\Dmt{\tilde{\Dm}}
\def\Se{\S^\e}
\def\se{\s^\e}
\def\ose{\overline{\se}}
\def\oS{\overline{\S}}
\def\tS{\tilde{\S}}
\def\tts{\hat{\S}}
\def\Sn{\S^\e}
\def\Sno{\overline{\Sn}}
\def\Snt{\tilde{\Sn}}
\def\Sm{\S^\eta}
\def\Smo{\overline{\Sm}}
\def\Smt{\tilde{\Sm}}
\def\ve{\vv^\e}
\def\ove{\overline{\ve}}
\def\vove{\ve\otimes\ve}
\def\vd{\vv^\delta}
\def\sd{\s^\delta}
\def\Sd{\S^\delta}
\def\Dd{\DD\vd}

\def\Wnd#1{W^{1,#1}_{\n, \diver}}
\def\Wndr{W^{1,r}_{\n, \diver}}

\def\o{\Omega}
\def\po{\partial \Omega}
\def\dt{\frac{\d}{\d t}}
\def\pt{\partial_t}
\def\it{\int_0^t \!}
\def\iT{\int_0^T \!}
\def\io{\int_{\o} \!}
\def\iq{\int_{Q} \!}
\def\iqt{\int_{Q^t} \!}
\def\ipo{\int_{\po} \!}
\def\ig{\int_{\Gamma} \!}
\def\igt{\int_{\Gamma^t} \!}

\def\d{\, \textrm{d}}

\def\mn{\mathcal{P}}
\def\du{\mathcal{W}}
\def\tr{\operatorname{tr}}
\def\tow{\rightharpoonup}

%------------------------------------------------

\newtheorem{theorem}{Theorem}[section]
\newtheorem{lemma}[theorem]{Lemma}
\newtheorem{proposition}[theorem]{Proposition}
\newtheorem{remark}[theorem]{Remark}
\newtheorem{corollary}[theorem]{Corollary}
\newtheorem{definition}[theorem]{Definition}
\newtheorem{example}[theorem]{Example}

\numberwithin{equation}{section}

\title[Implicitly constituted incompressible fluids]{On unsteady internal flows of incompressible fluids characterized by implicit constitutive equations in the bulk and on the boundary}

\thanks{M.~Bul\'i\v{c}ek and J.~M\'alek acknowledge the support of the project No. 20-11027X financed by the Czech Science foundation (GA\v{C}R). M.~Bul\'i\v{c}ek and J.~M\'alek are members of the {Ne\v{c}as} Center for Mathematical Modelling.}

\author[M.~Bul\'{i}\v{c}ek]{Miroslav Bul\'i\v{c}ek}
\address{Charles University, Faculty of Mathematics and Physics, Mathematical Institute, Sokolovsk\'{a} 83, 186~75, Prague, Czech Republic}
\email{mbul8060@karlin.mff.cuni.cz}

\author[J. M\'alek]{Josef M\'alek}
\address{Charles University, Faculty of Mathematics and Physics, Mathematical Institute, Sokolovsk\'{a} 83, 186~75, Prague, Czech Republic}
\email{malek@karlin.mff.cuni.cz}

\author[E.~Maringov\'{a}]{Erika Maringov\'{a}}
\address{Institute of Science and Technology Austria, Am Campus 1, 3400 Klosterneuburg, Austria}
\email{erika.maringova@ist.ac.at}


\begin{abstract}
Long-time and large-data existence of weak solutions for initial- and boundary-value problems concerning three-dimensional flows of \emph{incompressible} fluids is nowadays available not only for Navier--Stokes fluids but also for various fluid models where the relation between the Cauchy stress tensor and the symmetric part of the velocity gradient is \emph{nonlinear}. The majority of such studies however concerns models where such a dependence is \emph{explicit} (the stress is a function of the velocity gradient), which makes the class of studied models unduly restrictive. The same concerns boundary conditions, or more precisely the slipping mechanisms on the boundary, where the no-slip is still the most preferred condition considered in the literature. Our main objective is to develop a robust mathematical theory for unsteady internal flows of \emph{implicitly constituted} incompressible fluids with implicit relations between the tangential projections of the velocity and the normal traction on the boundary. The theory covers numerous rheological models used in chemistry, biorheology, polymer and food industry as well as in geomechanics. It also includes, as special cases, nonlinear slip as well as stick-slip boundary conditions. Unlike earlier studies, the conditions characterizing admissible classes of constitutive equations are expressed by means of tools of elementary calculus. In addition, a fully constructive proof (approximation scheme) is incorporated. Finally, we focus on the question of uniqueness of such weak solutions.
\end{abstract}

\maketitle

\section{Dedication} This article is written as a contribution to the celebration of the 100th anniversary of the birth of \emph{Olga Aleksandrovna Ladyzhenskaya} (March 7, 1922 - January 12, 2004) and to honor her scientific achievements.

Olga Ladyzhenskya seems to have been the first to initiate interest in the mathematical community to study incompressible fluid models that go beyond the Navier--Stokes equations. At the International Congress of Mathematicians in Moscow 1966, she presented arguments (see \cite{OAL67} and \cite{La72}), based on the kinetic formulation\footnote{The authors would be thankful for precise reference or notes confirming this.},
indicating that the macroscopic relation between the stress and the symmetric part of the velocity gradient should be polynomial. Ladyzhenskaya's model belongs to the class of power-law fluid models (sometimes also called modified or generalized Navier--Stokes fluids) characterized by a power-law index $r$, where the value $r=2$ corresponds to the Navier--Stokes fluid. Ladyzhenskaya was interested in the rigorous analysis of models with $r>2$: she has found that for those with $r\geq\tfrac52$ one could prove that the weak solution corresponding to the relevant initial- and boundary-value problem (in the sense of Leray and Hopf \cite{Leray1934, Hopf1951}) not only exists for long-time and large-data but it is unique. This uniqueness result should be contrasted with her counterexample to uniqueness of a weak solution of the Navier--Stokes equations in special time-dependent domains, see \cite{OAL_nonunique}. She also addressed, particularly in her subsequent studies, other aspects of weak solutions of these equations such as higher temporal and spatial differentiability and long-time behavior (the existence of a global attractor and estimates of its dimension).

K.~R.~Rajagopal together with the second author of this study reviewed Olga Ladyzhenskaya's foundational results concerning the analysis of modified Navier--Stokes equations\footnote{See also chapters in the book by J.-L. Lions \cite{Li69}.}, achieved during the period 1967--2003, in the second part of their handbook article~\cite{MaRa2005}. One of the objectives of this study is to provide a brief review of the results obtained in the mathematical analysis of fluids with nonlinear algebraic relation between the Cauchy stress and the velocity gradient achieved after 2003. The main objective is however to present a novel existence theory.

\section{Formulation of the problem and of the main result}\label{Sect2}

Materials are \emph{incompressible} if the volume of any measurable subpart of the body remains unchanged during a deformation process. For fluids flowing in a $d$-dimensional domain\footnote{Throughout the whole study, the term \emph{domain} stands for an open bounded connected set in $\mathbb{R}^d$.} $\Omega$, the condition of incompressibility expressed in the terms of the velocity $\vv = (v_1, \dots, v_d)$ takes the form
\begin{equation}
    \diver \vv(t,\x) =~ 0 \qquad \textrm{ for all } t \geq 0 \textrm{ and } \x\in \Omega.\label{pepa1}
\end{equation}
Incompressibility, which should be considered as a useful idealization, implies that the Cauchy stress tensor $\T$ is of the form
\begin{equation}
    \T = -p\I + \S, \label{pepa2}
\end{equation}
where only the part $\S$ can be determined experimentally. \emph{Homogeneous incompressible fluids} are characterized as incompressible fluids in which the density remains unchanged and is equal to a positive constant $\rho_*$. Such fluids automatically fulfil the balance of mass equation.

To conclude, setting  $Q=(0,T)\times \Omega$ and $\Gamma:= (0,T) \times \partial \Omega$ with $T>0$, the governing equations for unsteady flows of any homogeneous incompressible fluid flowing in a fixed domain $\Omega$ with no outflows and inflows and with initial velocity $\vv_0$ take the form
\begin{equation}
\label{pepa3}
\begin{aligned}
\diver \vv =&~ 0 &&\text{in } Q, \\
\rho_* (\pt \vv + \diver (\vv \otimes \vv) ) =&~ -\nabla p + \diver \S + \rho_* \b &&\text{in } Q, \\
\S=&~\S^T &&\text{in } Q,\\
\vv \cdot \n =&~0 &&\text{on } \Gamma, \\
\vv (0,\cdot) =&~\vv_0 &&\text{in } \o.
\end{aligned}
\end{equation}
Here $\b$ stands for the density of external body forces. The second equation in \eqref{pepa3} comes from the balance of linear momentum once \eqref{pepa2} is incorporated. The third equation says that the tensor $\S$ is symmetric; this implies that the balance of angular momentum is fulfilled. The fourth equation states that all considered flows are \emph{internal}: the fluid cannot enter or leave  $\Omega$.

The system \eqref{pepa3} is incomplete as, in $Q$, we have $d+1 + d(d+1)/2$ unknowns $\vv,p$ and $\S$, but merely $d+1$ equations. Also on $\Gamma$, we only have one scalar equation, but in fact $d$ boundary conditions are expected.

Taking the scalar product of the second equation in \eqref{pepa3} and $\vv$, integrating the result over $\Omega$ and using the remaining equations in \eqref{pepa3}, one obtains, after the integration over $(0,t)$ for any $t\in (0,T]$, the energy identity in the form (see \cite[Sect. 4.6]{MP2018} for details)
 \begin{equation}\label{pepa4}
    \frac{\rho_*}{2} \|\vv(t)\|_2^2 + \int_0^t\int_\Omega  \S:\D \d \x \d \tau + \int_0^t\int_{\partial \Omega} \s\cdot\vv_\tau \d S \d \tau = \int_0^{t} \int_\Omega \rho_*\b \cdot \vv  \d \x \d \tau + \frac{\rho_*}{2} \|\vv_0\|_2^2,
\end{equation}
where $\s$ stands for the projection of the normal traction $\T\n$ to the tangent plane, i.e.,
\begin{equation*}
\s:=-(\T\n)_\tau, \qquad \textrm{ where } \vec{z}_\tau := \vec{z} - (\vec{z} \cdot \n)\n,
\end{equation*}
$\n:\partial\Omega \to \mathbb{R}^d$ being the outer normal to $\partial\Omega$. Note that $(\T\n)_\tau
= (\S\n)_\tau$.

The second and third terms on the left-hand side represent two independent dissipation mechanisms: the former is associated with the internal friction inside the fluid, the latter corresponds to the interaction of the flowing fluid with (the inner part of) the boundary. Both terms should be, in accordance with the second law of thermodynamics, non-negative. For the Euler fluid, when $\S = \0$ and also $\s=\0$, both terms vanish. The first term will also vanish if all admissible flows are rigid (i.e. $\D = \0$), while the second term is equal to zero if all considered flows are subject to the no slip boundary condition, i.e.  $\vv_{\tau} = \0$ on $\Gamma$. For the Navier--Stokes fluid characterized by the constitutive equation
\begin{equation}
    \S = 2\nu_* \D, \qquad \textrm{ where } \nu_*>0, \label{pepa5}
\end{equation}
we conclude that\footnote{Analogously, for Navier's slip boundary condition
\begin{equation}
    \s = \gamma_* \vv_{\tau}, \qquad \textrm{ where } \gamma_*>0, \label{pepa7}
\end{equation}
we conclude that
\begin{equation}\label{pepa8}
    \int_{\partial \Omega} \s:\vv_{\tau} \d S = \frac{\gamma_*}{2} \int_{\partial \Omega} |\vv_{\tau}|^2 \d S + \frac{1}{2\gamma_*} \int_{\partial \Omega} |\s|^2 \d S.
\end{equation}}
\begin{equation}\label{pepa6}
    \int_\Omega \S:\D \d \x = \nu_* \int_\Omega |\D|^2 \d \x + \frac{1}{4\nu_*} \int_\Omega|\S|^2 \d \x.
\end{equation}
This in conjuction with the energy identity \eqref{pepa4} guarantees control of $\nabla \vv$ in $L^2(Q)$, which happens to be a key piece of information to establish long-time existence of a weak solution for any domain $\Omega$, $T>0$, $\nu_*>0$, $\b\in L^2(Q)$ and $\vv_0\in L^2(\Omega)$. This is what Leray and Hopf proved, see \cite{Leray1934} and \cite{Hopf1951}.

The above constitutive equations \eqref{pepa5} and \eqref{pepa7} are \emph{linear}. There are however many fluids (as is also illustrated in more detail in the next section) exhibiting \emph{nonlinear} relationship between $\S$ and $\D$. The same concerns slipping boundary conditions. Following Rajagopal \cite{Raj2003, Raj2006, RajSrin2008}, it is tempting to include all these equations under the umbrella of implicit equations relating $\S$ and $\D$ on the one hand and $\s$ and $\vv_{\tau}$ on the other hand. Hence, as nonlinear generalizations of \eqref{pepa5} and \eqref{pepa7}, we add to the problem \eqref{pepa3} the following equations
\begin{align}
\G(\S,\D) =&~\0 &&\text{in } Q, \label{pepa9}\\
\g(\s, \vv_\tau) =&~ \0 &&\text{on } \Gamma, \label{pepa10}
\end{align}
where $\G: \R^{d\times d} \times \R^{d\times d} \to \R^{d\times d}$ and $\g:\R^d \times \R^d \to \R^d$ are given continuous functions. Adding \eqref{pepa9} to the first three equations in \eqref{pepa3} we obtain a closed system of partial differential equations consisting of $d+1 + d(d+1)/2$ equations for  $d+1 + d(d+1)/2$ unknowns $\vv$, $p$ and $\S$. Adding \eqref{pepa10} to the fourth equation in \eqref{pepa3} we get $d$ equations on the boundary: one in the normal direction and $(d-1)$ of them are formulated at the tangent plane to $\partial \Omega$. Stated differently, the problem \eqref{pepa3} together with \eqref{pepa9} and \eqref{pepa10} is well-formulated. Motivated by Leray--Hopf's theory for the Navier--Stokes equations, it is natural to ask:
\smallskip
\begin{center}
{
\begin{minipage}{.8\linewidth}
    \emph{Can we formulate conditions on $\G$ and $\g$ that would allow us to  establish the long-time and large-data existence of a weak solution for \eqref{pepa3}, \eqref{pepa9} and \eqref{pepa10}? If so, can these conditions be formulated in terms of the tools of elementary calculus so that they are accessible to a broad scientific community?}
\end{minipage}
}
\end{center}
\smallskip

\noindent This study provides positive answers to these questions. Before formulating the main result, we give the admissibility conditions on $\G$ and $\g$. Here, we closely follow our preceding study~\cite{BMM21} focused however on a simpler problem (a ``mixed" formulation for problems of parabolic type).

Regarding the tensorial function $\G$ which determines the material response inside the domain $Q$, we assume that
\begin{itemize}
\item[(G1)] $\G$ is Lipschitz continuous, i.e. $\G \in \mathcal{C}^{0,1}(\R^{{d\times d}} \times \R^{{d\times d}})^{{d\times d}}$ and $\G(\0,\0)=\0$;
\item[(G2)] for almost all $(\S, \DD) \in \R^{{d\times d}} \times \R^{{d\times d}}$:
\begin{equation*}
\begin{split}
&\frac{\partial \G(\S,\DD)}{\partial \S}\geq 0, \quad \frac{\partial \G(\S,\DD)}{\partial \DD}\le 0, \quad \frac{\partial \G(\S,\DD)}{\partial \S}-\frac{\partial \G(\S,\DD)}{\partial \DD}>0, \\
&\textrm{and } \quad \frac{\partial \G(\S,\DD)}{\partial \DD}\left(\frac{\partial \G(\S,\DD)}{\partial \S}\right)^T\le 0;
\end{split}
\end{equation*}
\item[(G3)] one of the following holds:
\begin{equation*}
\begin{aligned}
&\text{either} &&\forall \DD\in \R^{{d\times d}} \quad \liminf_{|\S|\to +\infty} \G(\S, \DD) : \S > 0\\
&\text{or} &&\forall \S\in \R^{{d\times d}}\quad \limsup_{|\DD|\to +\infty} \G(\S, \DD) : \DD< 0;
\end{aligned}
\end{equation*}
\item[(G4)] there exist $C_1, C_2>0$ such that for all $(\S,\DD)\in \R^{{d\times d}} \times \R^{{d\times d}}$ fulfilling~$\G(\S,\DD)=\0$ we have
\begin{equation*}
\S :\DD \geq C_1 (|\S|^{r'} + |\DD|^r) -C_2, \qquad r':= {r}/({r-1}).
\end{equation*}
\end{itemize}
In (G2), we used the following notation for $(\G)_{ij}=G_{ij}$ and $(\S)_{ij}= S_{ij}$,
$$
\left(\frac{\partial \G}{\partial \S}\right)^{ij}_{kl}= \frac{\partial G_{ij}}{\partial S_{kl}}.
$$
Further, $\vec{A}^T$ denotes the transpose tensor to $\vec{A}$, i.e. $(\vec{A}^T)^{ij}_{kl}={A}_{ij}^{kl}$ and  $\vec{A}\vec{B}^T$ is the standard tensor multiplication, i.e. $(\vec{A}\vec{B}^T)^{ij}_{kl}=\sum_{m,n=1}^{d}{A}^{ij}_{mn}{B}^{kl}_{mn}$.  Also, for any tensor $\vec{A}\in \R^{d\times d}\times \R^{d\times d}$, the expression $\vec{A}\geq 0$ means that for any $\X\in \R^{d\times d}$ there holds
$$
\vec{A}\X : \X\geq 0 \qquad \textrm{(which, written in terms of components, is $\sum_{i,j,k,l=1}^d A^{ij}_{kl}X_{ij} X_{kl} \geq 0$)}.
$$
In addition, if we write $\vec{A}>0$ then we mean that the above inequality is strict for all $\X\neq \0$.

Before formulating similar conditions on $\g$, some comments are in order. First, note that the constitutive equation $\G(\S,\D)=\0$ can be replaced by $-\G(\S,\D)=\0$. Then all inequalities in (G2) and (G3) have the opposite signs except the last inequality in (G2). This ambiguity could be fixed for example by requiring that $\G$ is such that the first condition in (G2) holds (compare it with the special case of the Navier--Stokes fluids, see \eqref{pepa5}, when one would consider $\G(\S,\D)=\S - 2\nu_*\D$ and not $\G(\S,\D)= 2\nu_*\D - \S$). Second, as  the null points of $\G$ are of our interest, we can require the validity of (G2) only in the neighbourhood $\G(\S,\D) = \0$.
Also, the Lipschitz continuity in (G1) is required only to guarantee the existence of  the partial derivatives in (G2) almost everywhere. Alternatively, one can assume merely the continuity of $\G$ in (G1) and substitute (G2) by
\begin{equation}\tag{G2$^*$}\label{G2*}
    (\S_1 - \S_2):(\DD_1 - \DD_2) \geq 0 \quad \text{ for all } (\S_i, \DD_i) \text{ such that } \G(\S_i,\DD_i)=\0,~i=1,2.
\end{equation}
See also the statement of Lemma~\ref{pidilema} below. Finally, the conditions (G1)--(G4) are formulated using elementary tools of calculus (limes superior, partial derivatives). As proved in \cite{BMM21}, these conditions are equivalent to the statement that the set of null points of $\G$ is a maximal monotone $r$-coercive graph that passes through the origin (see \cite[Definition 3.1 and Lemma 3.2]{BMM21} for details.

Regarding the vectorial function $\g$ which determines the relation between the shear stress $\s$ and the tangential velocity $\vv_\tau$ on the boundary, we assume (in a similar way as above) that
\begin{itemize}
\item[(g1)] $\g \in \mathcal{C}^{0,1}(\R^d \times \R^d)^d$, and $\g(\0,\0)=\0$;
\item[(g2)] for almost all $(\s, \vv) \in \R^d \times \R^d$:
\begin{equation*}
\begin{split}
&\frac{\partial \g(\s,\vv)}{\partial \s}\geq 0, \quad \frac{\partial \g(\s,\vv)}{\partial \vv}\leq 0, \quad\frac{\partial \g(\s,\vv)}{\partial \s}-\frac{\partial \g(\s,\vv)}{\partial \vv}>0, \\
&\textrm{and } \quad \frac{\partial \g(\s,\vv)}{\partial \vv}\left(\frac{\partial \g(\s,\vv)}{\partial \s}\right)^T\leq 0;
\end{split}
\end{equation*}
\item[(g3)] one of the following conditions holds:
\begin{equation*}
\text{either} \quad\forall \vv\in \mathbb{R}^d \quad \liminf_{|\s|\to +\infty} \g(\s, \vv) \cdot \s > 0 \quad \text{or} \quad \forall \s\in \mathbb{R}^{d}\quad \limsup_{|\vv|\to +\infty} \g(\s, \vv) \cdot \vv< 0;
\end{equation*}
\item[(g4)] there exist $c_1, c_2>0$ such that, for all $(\s,\vv)\in \R^d \times \R^d$ fulfilling~$\g(\s,\vv)=\0$, the following condition holds:
\begin{equation*}
\s \cdot \vv \geq c_1 (|\s|^{q'} + |\vv|^q) -c_2, \qquad q':= {q}/({q-1}).
\end{equation*}
\end{itemize}
The above comments related to $\G$ are applicable to $\g$ as well.
%In (g1)--(g4), we used the notation which is analogous to the one in (G1)--(G4). Also, since we can replace $\g$ by $-\g$, it is clear that all inequalities in (g2) and (g3) can be equivalently formulated with the opposite sign except the last inequality in (g2). Also, note that the following holds\footnote{See Lemma~\ref{pidilema} for details.}
%\begin{equation}\tag{g2$^*$}\label{g2*}
%    (\s_1 - \s_2):(\vv_1 - \vv_2) \geq 0 \quad \text{ for all } (\s_i, \vv_i) \text{ such that } \g(\s_i,\vv_i)=\0,~i=1,2,
%\end{equation}
%which, together with only continuity of $\g$ (instead of assumed Lipschitz continuity) could be considered as an alternative to (g2).

Now, we are ready to formulate our main result (in a vague way):
\smallskip
\begin{center}
\fbox
{
\begin{minipage}{.8\linewidth}
    \emph{For arbitrary $\o$, $T$, $\vv_0$ and $\b$, and for any $\G$ and $\g$ fulfilling (G1)--(G4) with $r>\frac{2d}{d+2}$ and (g1)--(g4) with $q>1$, there exists a weak solution %$\vv$, $\S$, $\s$ and $p$ that solve
  to the problem \eqref{pepa3}, \eqref{pepa9} and \eqref{pepa10}.}
\end{minipage}
}
\end{center}
\bigskip

{\bf Overview of the existence theory for $\G(\S,\D)=\0$ after 2003}.  We restrict our discussion to the most interesting case $d=3$. In the first period, prior 2011, the focus of research were models of the type $\S = (1 + |\D|^2)^{\frac{r-2}{2}}\D$ following the goal to establish global-in-time existence theory for large data for models with low values of $r$. (Note that the local-in-time existence of smooth solutions to models of power-law type is addressed in \cite{BoPr2007}.) Let us recall that the problems studied by Olga Ladyzhenskaya concerned the subcritical regime when the velocity itself is an admissible test function in the weak formulation of the balance of linear momentum. This corresponds to the case when $r\geq \frac{11}{5}$. The method of Lipschitz truncation developed in \cite{FMS03} (see also \cite{DMS09}) for time-independent (stationary) problems covers the case $r>\frac{6}{5}$ but was left open for the evolutionary case. For the evolutionary case, the ``best" result known around the year 2005, see \cite{FMS00}, covered the case $r\geq \frac{8}{5}$ using the $L^{\infty}$-truncation technique; the approach is restricted to the spatially periodic problem. The extension for flows in general bounded domains subject to Navier's slip boundary was established in \cite{BuMaRa07}, while the no-slip boundary conditions were successfully treated in \cite{wolf}, still for $r\geq \frac{8}{5}$. Finally, Diening, R\accent23u\v{z}i\v{c}ka and Wolf, see \cite{DiRuWo10}, inspired by the works of Kinnunen and Lewis \cite{LeKin02}, extend the method of Lipschitz truncation to the evolutionary case and proved the existence of a weak solution to the evolutionary problem with no-slip boundary condition for $r> \frac{6}{5}$. The remaining case $r\in [1, \frac{2d}{d+2}]$ is covered by two somehow contradictory recent results, see \cite{AF} and \cite{MR4328053}, which can be interpreted in the way that the range of possible $r$'s in Theorem~\ref{genresult} is optimal. In fact, Abbatiello and Feireisl \cite{AF} introduce a novel generalized concept of solution (dissipative solution) and establish its existence theory for $r\in (1, \frac{2d}{d+2}]$. The theory is developed for a smaller class of possible constitutive relations than considered here. More importantly, their concept of solution does not imply either the validity of weak formulation of balance of linear momentum  (see \eqref{WF} in Sect.~\ref{Sect5}) or the validity of $\G(\S,\D) = \0$ almost everywhere in $Q$. On the other hand, in \cite{MR4328053}, Burczak, Modena and Sz\'{e}kelyhidi show the non-uniqueness of even Leray--Hopf solutions for $r<\frac{2d}{d+2}$. From this perspective, it seems also to be reasonable to consider only the case $r\geq \frac{2d}{d+2}$. Note that even for this range of $r$'s Burczak et al. \cite{MR4328053} prove the result concerning non-uniqueness of a very weak\footnote{Very weak solutions are those that do not belong to the natural energy function space.} solution. As we are dealing with weak solutions, their result is not applicable to our setting.

Inspired by the foundational works on implicit constitutive relations, see \cite{Raj2003,Raj2006}, the question to develop a robust theory covering the whole class of implicitly constituted incompressible fluids arose. Following initial attempts (see \cite{MalekETNA2008, BGMS1}), a successful theory covering both polynomial and activated fluids was established in \cite{BGMS2, BGMS3}, even in a broader context than considered here: the $r$-coercivity condition is generalized in terms of Young's functions in the setting of Orlicz spaces and the constitutive equation \eqref{pepa9} was allowed to vary with time and space, i.e. $\G(t,\x,\S,\D) = \0$ in $Q$. As discussed in length in \cite[pp 2048--2049]{BMM21}, there are two shortcomings of the results proved in \cite{BGMS2, BGMS3} (a non-constructive proof and an a~priori assumption concerning the existence of a~Borel measurable selection). These shortcomings motivated the development of an alternative approach, see \cite{BMM21}. The extension of the approach developed in \cite{BMM21} for problems of parabolic type to problems involving flows of incompressible fluids is one of the main objections of this study.

Note that the assumption (g4) eliminates no-slip and perfect-slip boundary conditions from the analysis presented here. However, as indicated in the above discussion of available results one can incorporate both conditions into the analysis. The case of $\s = \0$ on $\Gamma$ is in fact easy as the boundary term just vanishes. For no-slip boundary conditions, one needs to change the function space for the velocity and pay attention to differences associated with the reconstruction of the pressure (see \cite{wolf} and \cite{BlMaRa2020} for details).

Numerical analysis of finite-element based discretizations of problems of the type \eqref{pepa3}, \eqref{pepa9} and \eqref{pepa10}, completed with computational experiments, is addressed in \cite{Diening2013, Kreuzer2016, Farrell2020, Farrell2022, Suli2020, Heid2022}.

{\bf Structure of the paper.} In Sect.~\ref{Sect3}, we illustrate how rich the classes of fluids under consideration are by providing a list of models used in various areas of science (completed by a list of references). In Sect.~\ref{Sect4}, we introduce, in a constructive way, $\e$-approximations of the constitutive equations and provide a summary of their properties (proved in~\cite{BMM21}). In particular, at this approximate level the term $-\diver\S$ leads to Lipschitz continuous uniformly monotone elliptic operator (the nicest one can wish to deal with). After introducing basic function spaces in Section~\ref{Sect5} we give a  precise formulation of the main theorem including also the precise definition of weak solution to \eqref{pepa3}, \eqref{pepa9} and \eqref{pepa10}. Here we also recall properties of the Lipschitz approximations of Bochner functions needed in the proof of the main theorem. This forms the content of Sect.~\ref{Sect6}. As this article aims at surveying the results in the field, we give, in Sect.~\ref{Sect7}, a summary of the results that concern similar problems including an additional component that makes the whole problem more complicated. Finally, we comment on the available uniqueness results regarding the studied problem in Sect.~\ref{Sect8}.

\section{Examples of implicit constitutive equations} %$\G(\S,\D)=\0$ and $\g(\s,\vv_{\tau})=\0$}
\label{Sect3}

The purpose of this section is to provide an illustrative list of models and boundary conditions covered by the implicit equations \eqref{pepa9} and \eqref{pepa10}.
The aim is to show that these classes of fluids and boundary conditions are rich and particular models appear in various areas of science and engineering. We first focus on the constitutive equations in the bulk, then we discuss the boundary conditions.

{\bf Constitutive equations covered by \eqref{pepa9}.} The fact that various fluids exhibit a nonlinear rheological relation between the shear stress and the shear rate was known already at the end of the 19th century, see Schwedoff \cite{Schwedoff1890}, Troutan \cite{trouton1906} and further references in books on non-Newtonian fluids, such as Bird, Amstrong, Hassager
\cite{bird77}, Huilgol \cite{huilgol}, Schowalter
\cite{schowalter}, or in the survey paper \cite{MRR95}. There are hundreds of models belonging to this framework used in chemistry, biofluid rheology, geomechanics, food industry, etc. A discussion of various aspects of these models can be found in~\cite{MRR95}, with references to the chemical engineering and colloidal mechanics literature (\cite{Carreau1972}, \cite{CHRISTIANSEN1973}, \cite{Powell1944}, \cite{Ree1958}, \cite{Sutterby1965}, \cite{TURIAN1969}), ice-mechanics and glaciology (\cite{Kjartanson1988}, \cite{METZNER1956}, \cite{whillans1993}), blood rheology (\cite{Cho1989}, \cite{Cho1991}, \cite{Cokelet1972}, \cite{Cross1965}, \cite{davies1990}, \cite{Nakamura1988}, \cite{Powell1944}, \cite{Steffan1990}, \cite{Walawender1975}, \cite{casson1959, cebral2005, Quemada1978, YELESWARAPU1998}, \cite{Galdi2008} and \cite{Fasano2017}).

These models fit to the setting characterized by the form
\begin{equation}\label{pepa11}
\T = -p\I + \S, \quad \text{ where } \quad \alpha(|\S|^2, |\D|^2)\S = 2 \nu(|\S|^2, |\D|^2)\D.
\end{equation}
As the fluid is incompressible, and consequently the trace of $\D$ vanishes, one observes that within the class \eqref{pepa11} one has $p = - \frac13\operatorname{tr}\T$.

\begin{table}[ht]
\begin{tabular}{|l|l|l|}
\hline
\rule{0in}{.15in}Model & $\nu(|\D|^2)$ & $\nu(|\S|^2)$
\rule[-.1in]{0in}{0in}\\ \hline\hline
\rule{0in}{.15in}Ostwald-de Waele \cite{Ostwald1925},\cite{waele1923}  & $\nu_0 |\D|^{m-1} $ & \\
Glen  \cite{Glen1955} & & $A |\S|^{m-1}$ \\
Carreau \cite{Carreau1972} \quad & $\nu_\infty + \frac{\nu_0 - \nu_\infty}{(1+A|\D|^2)^{n/2}}$ & \\
Blatter \cite{Pettit2003}, \cite{blatter1995} & & $\frac{A}{(|\S|^2 + \tau_0^2)^{(n-1)/2}}$ \\
Carreau-Yasuda \cite{yasuda1979}
\quad & $\nu_\infty + \frac{\nu_0 - \nu_\infty}{(1+A|\D|^a)^{n/a}}$ & \\
Eyring \cite{Eyring1936} & $\nu_\infty + (\nu_0 - \nu_\infty) \frac{\operatorname{arcsinh} (A|\D|)}{A|\D|}$ & \\
Sisko \cite{Sisko1958} & $\nu_\infty + A|\D|^{n-1}$ &\\
Cross \cite{Cross1965} & $\nu_\infty + \frac{\nu_0 - \nu_\infty}{1+A|\D|^n}$ & \\
Ellis \cite{Matsuhisa1965} & & $\frac{\nu_0}{1+A|\S|^{n-1}}$ \\
Seely \cite{Seely1964} & & $\nu_\infty + (\nu_0 - \nu_\infty) e^{-|\S|/\tau^2_0}$
 \rule[-.1in]{0in}{0in} \\ \hline
\end{tabular}\vspace{.3cm}
\caption{Frequently used models in material sciences, chemical engineering, biomechanics and geophysics. Here $\nu_0$, $\nu_\infty$, $m$, $a$ and $A$ are positive constants, while $n$ and $\tau_0$ are real numbers. The models are taken from \cite[Section 4.5]{MaPr2018}.}
\label{tab:chem}
\end{table}

In Table~\ref{tab:chem}, we distinguish two special subclasses of \eqref{pepa11}, namely $\S = 2\nu(|\D|^2)\D$ and $\S = 2\nu(|\S|^2)\D$.
The simplest deviation from the Navier--Stokes fluid model represents the power law model that can be described in two equivalent ways as follows:
\begin{equation}
    \label{pepa12}
    \S = 2\nu_0 |\D|^{r-2}\D \quad\iff \quad \D = \frac{1}{(2\nu_0)^\frac{1}{r-1}} |\S|^{\frac{2-r}{r-1}} \S\,,
\end{equation}
where $r\in (1,\infty)$ and $\nu_0>0$. Referring to Table~\ref{tab:chem}, we thus observe that the same model is called Ostwald-de Waele's model in chemistry, while it is named Glen's model in geomechanics. Denoting $r':=r/(r-1)$ we also have
\begin{equation}
    \label{pepa13}
    \S:\D = \left(\frac{1}{r} + \frac{1}{r'}\right) \S:\D = \frac{2\nu_0}{r} |\D|^r + \frac{r-1}{r(2\nu_0)^{\frac{1}{r-1}}}|\S|^{r'}\,,
\end{equation}
which reduces to \eqref{pepa6} if $r=2$. It also serves as the main motivation for the $(r,r')$-coercivity assumption (G4). All the models listed in Table~\ref{tab:chem} describe, for suitable range of parameters, a non-Newtonian phenomenon called \emph{shear thinning/shear thickening} (the generalized viscosity is decreasing/increasing function of the shear rate).

The constitutive equations of the form \eqref{pepa9} are also suitable to describe fluids with the activation criteria. Bingham and Herschel--Bulkey fluids \cite{bingham1922,Herschel1926} can be written in the form
\begin{equation}\label{pepa15}
\D = \frac{1}{2\nu(|\D|^2)} \frac{(|\S| - \tau_*)^{+}}{|\S|} \S, \qquad \tau_*\in (0,\infty),
\end{equation}
where $\nu$ is a positive constant in the case of Bingham fluids and is a polynomial function of $\D$ in the case of Herschel--Bulkley fluids. It is proved in \cite[Appendix, Example A.3]{BMM21} that Bingham fluids satisfy (G1)--(G4) with $r=2$. Following the same line of argument, one can check that Herschel--Bulkley fluids with $v(|\D|^2) = (1 + |\D|^2)^{\frac{r-2}{2}}$ also satisfy (G1)--(G4).

Activated Euler fluids, see \cite{BlMaRa2020}, are described by the formula
\begin{equation}\label{pepa16}
\S = 2\nu(|\D|^2) \frac{(|\D| - \delta_*)^{+}}{|\D|} \D, \qquad \delta_*\in (0,\infty).
\end{equation}
If $\nu$ is constant, then the fluid behaves as the Navier--Stokes fluid once $|\D|$ exceeds the activation parameter $\delta_*$. It is straightforward to check that these models fulfil (G1)--(G4). (A large-data analysis of activated Euler fluids, for steady and unsteady flows and for various boundary conditions including complete slip as well as no-slip is developed in \cite{BlMaRa2020}.)

The examples discussed above are summed up in the following Table~\ref{tab:1}, where for brevity we set all physical constants to be $1$, except the exponent $r$ related to the $(r,r')$-coercivity condition (G4). If $r$ does not appear in the equation, then the model leads to (G4) with $r=2$. %(compare with Table 1 in~\cite{BMM21}).
\begin{table}[h]\label{table2}
\begin{tabular}{|l@{\hspace{.2in}}|l@{\hspace{.2in}}|}
\hline
\rule{0in}{.15in}$\S = \S^*(\D)$ & $\D = \DD^*(\S)$
\rule[-.1in]{0in}{0in}\\ \hline\hline
\rule{0in}{.15in}$\S = \D$ & $\D = \S$ \\
$\S = |\D|^{r-2} \D$ & $\D = |\S|^{r'-2} \S$ \\ %\hline
$\S = (1+ |\D|)^{r-2} \D$ &  $\D = (1+ |\S|)^{r'-2}\S$ \\
$\S = (1+ |\D|^2)^{\frac{r-2}{2}} \D$ &  $\D = (1+ |\S|^2)^{\frac{r'-2}{2}}\S$ \\
$\S = (|\D|-\delta_*)^+\frac{\D}{|\D|}$ \qquad $\delta_*>0$ & $\D = (|\S|-\tau_*)^+ \frac{\S}{|\S|}$ \qquad $\tau_*>0$
\rule[-.1in]{0in}{0in} \\ \hline
\end{tabular}\vspace{.3cm}
\caption{Examples of two classes of explicit constitutive relations covered by \eqref{pepa9}. The first two lines describe the Navier--Stokes and the power-law fluids with the power-law index $r\in (1,\infty)$, $r':=r/(r-1)$. In these rows, the formulas in the two columns are equivalent, see \eqref{pepa12}. The formulas in the third and fourth row hold for $r\in (-\infty,\infty)$; in the range $r>1$ the models in the left and right columns behave in the same way for large values of $|\D|$ and $|\S|$; for $r\geq 1$ and $r'\geq 1$ these formulas satisfy (G2) or (G2*), i.e. the response is monotone. For $r<1$ and for $r'<1$ the formulas in the left and right column behave differently; their response is non-monotone, see \cite{mpr10, LerouxKRR, PerPr, Janecka2019} for details. The last row describes an activated Euler fluid (left) and a Bingham fluid (right); after the activation takes place, both fluids behave as a Navier-Stokes fluid.}
\label{tab:1}
\end{table}

In conclusion, constitutive equation \eqref{pepa9} covers models designed to describe two non-Newtonian phenomena:  \emph{shear thinning/shear thickening} and \emph{the presence of activation criteria in a simple shear flow}. (A detailed description of non-Newtonian phenomena is given for example in \cite{MaRa2005}.) Interestingly, \eqref{pepa9}  covers one additional phenomenon, called \emph{normal stress differences}, which is usually attributed to viscoelastic nature of the fluid; see Perl\'{a}cov\'{a} and Pr\accent23u\v{s}a \cite{PerPr} for more details.
\bigskip

{\bf Constitutive equations (boundary conditions) covered by \eqref{pepa10}.} Boundary equations are of the same importance as the constitutive equations in the bulk. This assertion can be supported by a recent study \cite{Chab2022}, where flows are shown to change quantitatively in an essential manner in dependence of the boundary conditions (only linear Navier's slip boundary conditions were tested with their limiting cases (no slip vs complete slip).

Navier proposed a linear constitutive relation \eqref{pepa5} as the proper boundary condition in \cite{Navier1823}. Stokes \cite{Stokes1845} discusses the boundary conditions at length and one variant that he considers concerns a nearly quadratic relation between wall shear stress and the velocity. He states: ``\textit{\ldots when the velocity is not small the tangential force called into action by the sliding of water over the inner surface of the pipe varies nearly as the square of the velocity}". Mooney \cite{Mooney1931} proposed a more general form of slip and introduced a technique that evaluates this relationship. Comprehensive overviews concerning general boundary conditions and slipping mechanisms can be found for example in \cite{Hatzi} or \cite{RaRa99}. More detailed discussions concerning the boundary conditions, their importance, including references to earlier studies are available in \cite[A.4]{MaRa2005} and \cite[Sect. 4.6]{MaPr2018}. An overview of basic models of the type \eqref{pepa10} is given in Table~\ref{tab:bc}.

\begin{table}[h]
\begin{tabular}{|l@{\hspace{.2in}}|l@{\hspace{.2in}}|}
\hline
\rule{0in}{.15in}$\s = \s^*(\vv_\tau)$ & $\vv_\tau = \vec{d}^*(\s)$
\rule[-.1in]{0in}{0in}\\ \hline\hline
$\s = \vv_\tau$ & $\vv_\tau = \s$ \\
$\s = |\vv_\tau|^{q-2} \vv_\tau$ & $\vv_\tau = |\s|^{q'-2} \s$ \\ %\hline
$\s = (1+ |\vv_\tau|)^{q-2} \vv_\tau$ &  $\vv_\tau = (1+ |\s|)^{q'-2}\s$ \\
$\s = (1+ |\vv_\tau|^2)^{\frac{r-2}{2}} \D$ &  $\vv_\tau = (1+ |\s|^2)^{\frac{q'-2}{2}}\s$ \\
$\s = (|\vv_\tau|-\beta_*)^+\frac{\vv_\tau}{|\vv_\tau|} \qquad \beta_*>0$ & $\vv_\tau = (|\s|-\sigma_*)^+ \frac{\s}{|\s|} \qquad \sigma_*>0$
\rule[-.1in]{0in}{0in}\\ \hline
\end{tabular}\vspace{.3cm}
\caption{Examples of two classes of boundary conditions belonging to \eqref{pepa10}. The first line describes Navier's slip. The second, third and fourth rows describe nonlinear slip of polynomial type with the power-law index $q\in (1,\infty)$, $q':=q/(q-1)$. The last row describes the stick-slip boundary condition (right column) and the boundary condition that describes the complete slip before activation and Navier's slip once the activation takes place (left column).}
\label{tab:bc}
\end{table}

Measurements for molten polymers clearly document that there are nonlinear responses between $\s$ and $\vv_\tau$ including various activations. For example, Hatzikiriakos in \cite{Hatzi} considers models of power-law type, i.e.,
\begin{equation}\label{powerBC}
\vv_\tau =  \gamma |\s|^{q'-2}\s,
\end{equation}
referring to \cite{Ramamurthy1986} for the value $q'=3$, to \cite{Hatzi1991} for $q'=4$ and to \cite{Hill1990}, \cite{Kalika1987} for $q'=7$. Nonlinear responses \eqref{powerBC} include the material parameter $\gamma$ that can be a function of other relevant quantities. Besides \cite{Hatzi}, models of the type \eqref{powerBC} were studied also in \cite{Hatzi1992}, \cite{Chauff1979}, \cite{Lau1986}, \cite{Cohen1985}.

Stick-slip boundary conditions were added to large-data and long-time existence analysis in \cite{BuMa16, BuMa2017}. A treatment of complex nonlinear (non-monotone) boundary conditions of stick-slip type is presented within the context of analysis of Kolmogorov's two-equation model of turbulence in \cite{BuMa19}. How the choice of boundary conditions influences the definition of proper function spaces and the subsequent analysis is studied for different boundary conditions within the context of activated Euler's fluids in \cite{BlMaRa2020}.

%- Hatzikiriakos \cite{Hatzi} power-law slip model with references, even suggestion of power-law dynamic slip model - equation (19); some references coincide with \cite{RaRa99} \\
%- Journal of Rheology, Journal of Fluid Mechanics, Journal of nonNewt Fluid Mechanics\\
%-zpet ke korenum \cite{RaRa99} tam je v apendixu spostu modelu s nelinearnim slipem a odkazy nafyzikalni literatury. \\

\section{\texorpdfstring{$\e$}{e}-approximations of implicit constitutive equations}\label{Sect4}%$\G(\S,\D)=\0$ and $\g(\s,\vv_{\tau})=\0$}

%First, we define the maximal monotone graph.
%\begin{definition}[Maximal monotone $r$-coercive graph]\label{maxmongrafA}
%Let $r \in (1, \infty)$ and $r' := \frac{r}{r-1}$. We say that a subset $\A$ of  %$\R^{{d\times d}} \times \R^{{d\times d}}$ is a maximal monotone $r$-coercive graph %if
%\begin{itemize}
%\item[(A1)] $(\0,\0) \in \A$;
%\item[(A2)] For any $(\S_1,\DD_1), (\S_2,\DD_2) \in \A$
%\begin{equation*}
%(\S_1 - \S_2):(\DD_1 - \DD_2) \geq 0;
%\end{equation*}
%\item[(A3)] If for some $(\S,\DD) \in \R^{{d\times d}} \times \R^{{d\times d}}$ and for all $(\oS,\oD) \in \A$
%begin{equation*}
%(\S - \oS):(\DD - \oD) \geq 0,
%\end{equation*}
%then $(\S,\DD) \in \A$;
%\item[(A4)] There exist $C_1, C_2 >0$ such that for all $(\S,\DD) \in \A$
%\begin{equation*}
%\S :\DD \geq C_1(|\S|^{r'} +|\DD|^{r})-C_2.
%\end{equation*}
%\end{itemize}
%\end{definition}
%The condition (A1) means that $\A$ passes through the origin, (A2) states that the graph $\A$ is monotone, while (A3) states that $\A$ is maximal monotone. Finally, (A4) states that the graph $\A$ is $r$-coercive.

In our preceding study~\cite{BMM21}, while studying problems of parabolic type (think of a nonlinear heat equation), we found out the structural assumptions on the implicit function that characterizes the relation between the (heat) flux and the (temperature) gradient which allows us to built a theory parallel (``equivalent" - in the sense specified in \cite{BMM21}) to that of the maximal monotone $r$-coercive (potentially multi-valued) graphs. Here, we intentionally prefer to avoid using the concept of maximal monotone graphs and we wish to  present a theory based only on the assumptions on $\G$ and $\g$ that require only the knowledge of basic tools of calculus.

The intention of this section is to introduce $\e$-approximations of the functions $\G$ and $\g$ and summarize their nice properties: the approximations always lead to $L^2$-coercivity of $\S$ and $\D$ (and similarly for $\s$ and $\vv_{\tau}$), on the $\e$-approximation level, $\S$ is always a function of $\D$ and upon inserting it into $-\diver\S$ one obtains a Lipschitz continuous uniformly monotone operator.   Also here, we follow closely the approach developed in Sect. 4 in~\cite{BMM21}.

\begin{lemma}%[compare with Section 4 in~\cite{BMM21}]
    \label{convGg} Let $\G$ satisfy (G1)--(G4) for any $r>1$, let $\g$ satisfy (g1)--(g4) for any $q>1$ and let $\e\in(0,1)$. Then the approximating functions defined by
    \begin{subequations}\label{def:Gege}
    \begin{align}
        \Ge(\S,\DD) &:= \G(\S-\e\DD, \DD-\e\S), \label{def:Ge}\\
        \ge(\s,\vv) &:= \g(\s-\e\vv, \vv-\e\s) \label{def:ge}
    \end{align}
    \end{subequations}
    satisfy (G1)--(G4) and (g1)--(g4) with $r=q=2$. Also, there exist $\tilde{C}_1, \tilde{C}_2, \tilde{c}_1, \tilde{c}_2 >0$ independent of $\e$ such that
    \begin{subequations}
    \begin{align}
        \Ge(\Se,\DDe) = \0  \quad &\implies \quad ~\Se :\DDe \geq \tilde{C}_1 (|\Se|^{\min\{2,r'\}} + |\DDe|^{\min\{2,r\}}) - \tilde{C}_2, \label{indeG}\\
        \ge(\se,\ve) = \0  \quad &\implies \quad  \se \cdot \ve \geq \tilde{c}_1 (|\se|^{\min\{2,q'\}} + |\ve|^{\min\{2,q\}}) - \tilde{c}_2.\label{indeg}
    \end{align}
    \end{subequations}
    Moreover, there exist two unique functions (single-valued mappings)
    \begin{equation}\label{selectionSs}
        \S^*_{\!\e} : \R^{d\times d}\to \R^{d\times d} \qquad \text{ and } \qquad \s^*_\e : \R^d \to \R^d
    \end{equation}
    satisfying
%        \begin{align*}
%            \Ge(\S, \DD) = \0 &\iff \S = \S^*_{\!\e} (\DD), \\
%            \ge(\s, \vv) = \0 &\iff \s = \s^*_\e (\vv), \\
%            \S^*_{\!\e}(\0)=\0 \qquad &\text{ and } \qquad \s^*_\e(\0)=\0.
%        \end{align*}
        \begin{align*}
            \Ge(\S, \DD) = \0 \iff \S = \S^*_{\!\e} (\DD), \qquad
            \ge(\s, \vv) = \0 \iff \s = \s^*_\e (\vv), \qquad
            \S^*_{\!\e}(\0)=\0, \quad \s^*_\e(\0)=\0,
        \end{align*}
    and both $\S^*_{\!\e}$, $\s^*_\e$ are Lipschitz continuous and uniformly monotone, i.e. there exist positive constants $C_1(\e)$, $C_2(\e)$, $c_1(\e)$, $c_2(\e)>0$ such that, for any $\DD_1, \DD_2 \in \R^{d\times d}$ and any $\vv_1, \vv_2 \in \R^d$,
    \begin{subequations}
    \begin{align}
        |\S^*_{\!\e}(\DD_1) - \S^*_{\!\e}(\DD_2)| &\leq C_1(\e)|\DD_1-\DD_2|, \\
        (\S^*_{\!\e}(\DD_1) - \S^*_{\!\e}(\DD_2)):(\DD_1-\DD_2) &\geq C_2(\e)|\DD_1-\DD_2|^2, \\
        |\s^*_\e(\vv_1) - \s^*_\e(\vv_2)| &\leq c_1(\e)|\vv_1-\vv_2|, \\
        (\s^*_\e(\vv_1) - \s^*_\e(\vv_2))\cdot(\vv_1-\vv_2) &\geq c_2(\e)|\vv_1-\vv_2|^2. \label{nonnegge}
    \end{align}
    \end{subequations}
    If, in addition, for any bounded measurable $U\subset Q$ and for $\Se, \DDe: U\to \R^{d \times d}$
    \begin{equation}\label{bddSD}
        \Ge(\Se,\DDe) = \0  \text{ a.e. in } U \quad \text{ and } \quad \limsup_{\e\to 0_+} \int_U \Se:\DDe \d \x \d t \leq C
    \end{equation}
    then there exist $\S \in L^{r'}(U)$ and $\DD \in L^r(U)$ so that (for subsequences)
    \begin{align*}
        \Se &\tow \S &&\text{weakly in } L^{\min\{2,r'\}}(U) \\
        \DDe &\tow \DD &&\text{weakly in } L^{\min\{2,r\}}(U).
    \end{align*}
    Moreover, if
    \begin{equation}\label{bddSeDe}
        \limsup_{\e\to 0_+} \int_U \Se:\DDe \d \x \d t \leq \int_U \S:\DD \d \x \d t,
    \end{equation}
    then
    \begin{equation*}
        \G(\S,\DD) = \0 \text{ a.e. in } U \quad \text{ and } \quad \Se:\DDe \tow \S:\DD \text{ weakly in } L^1(U).
    \end{equation*}
    Analogously, if for a bounded and measurable $V\subset \Gamma$ and for $\se,\ve: V \to \R^d$
    \begin{equation}\label{bddsv}
        \ge(\se,\ve) = \0  \text{ a.e. in } V \quad \text{ and } \quad \limsup_{\e\to 0_+} \int_V \se\cdot \ve \d S \d t \leq C,
    \end{equation}
    then there exist $\s\in L^{q'}(V)$ and $\vv \in L^q(V)$ so that (for subsequences)
    \begin{align*}
        \se &\tow \s &&\text{weakly in } L^{\min\{2,q'\}}(V), \\
        \ve &\tow \vv &&\text{weakly in } L^{\min\{2,q\}}(V);
    \end{align*}
    moreover, if
    \begin{equation}\label{bddseve}
        \limsup_{\e\to 0_+} \int_V \se\cdot \ve \d S \d t \leq \int_V \s\cdot\vv \d S \d t,
    \end{equation}
    then
    \begin{equation*}
        \g(\s,\vv) = \0 \text{ a.e. in } V \quad \text{ and } \quad \se\cdot \ve \tow \s \cdot \vv \text{ weakly in } L^1(V).
    \end{equation*}
\end{lemma}
\begin{proof}
    See~\cite{BMM21}, Lemma~4.1, Lemma~4.2 and Lemma~4.4.
\end{proof}

\begin{lemma}\label{pidilema}
Let $\G$ satisfy (G1)--(G4) with $r\in(1,\infty)$. For every $\DD \in L^r(Q)$ there exists an $\S \in L^{r'}(Q)$ such that $\G(\S, \DD)=\0$ a.e. in $Q$. Moreover,~\eqref{G2*} holds.
\end{lemma}
\begin{proof}
See~\cite{BMM21}, Lemma~4.5 and Lemma~4.3. Property (G2$^*$) follows from the fact that the null points of $\G$ generate a maximal monotone graph, see~\cite[Lemma~3.2]{BMM21}.
\end{proof}

\section{Notation, function spaces and precise formulation of the main result}\label{Sect5}

\textbf{Notation}. Let $\Omega\subset\mathbb{R}^d$ be a domain. We say that $\Omega$ is a Lipschitz domain/$\C^{1,1}$-domain and we write $\Omega\in \C^{0,1}$/$\Omega\in \C^{1,1}$ if, roughly speaking, the boundary $\partial\Omega$ can be covered by finite number of overlapping $\C^{0,1}$/$\C^{1,1}$ mappings. For $t \in (0,T]$, we denote $Q^t:=[0,t) \times \o$ and $\Gamma^t:=[0,t) \times \Gamma$, and we recall that $Q := Q^T$ and $\Gamma := \Gamma^T$. The abbreviation \emph{a.a.} stands for \emph{almost all}, while \emph{a.e.} stands for \emph{almost everywhere}. Generic constants, that depend only on the data but are independent of any approximation parameter, are denoted by~$C$ and may vary from line to line.

For a Banach space~$(X,\|\!\cdot\!\|_X)$, its dual is denoted by $X^*$. For $x \in X$ and $x^*\in X^*$, the duality is denoted by $\langle x^*, x \rangle_X$. For $r \in [1,\infty]$, we denote $(L^r(\o), \|\!\cdot\!\|_{r})$ and $(W^{1,r}(\o),\|\!\cdot\!\|_{W^{1,r}(\Omega)})$ the corresponding Lebesgue and Sobolev spaces with the norms defined in the standard way.
%\begin{align*}
%\| f \|_{L^r(\o)} &:=
%\begin{cases}
%\left( \io |f|^r \d \x \right)^{\frac{1}{r}} &\text{if } r \in [1, \infty), \\
%\esssup_{x \in \o} |f(x)| &\text{if } r = \infty,
%\end{cases} \\
%\hspace{2cm}\| f \|_{W^{1,r}(\o)} &:= \|f\|_{L^r(\o)} +  \|\nabla f \|_{L^r(\o)}.
%\end{align*}
Bochner spaces are denoted by $L^r(0,T;X)$. We use the notation $L^r(\o;\R^d)$ and $L^r(\o;\R^{{d\times d}})$ for Lebesgue spaces of vector- or matrix-valued functions, respectively. $\C^{\infty}_0(U)$ stands for smooth functions with compact support in an open set $U$.

Next, we define the function spaces of divergenceless functions with the normal component vanishing on the boundary that are relevant to our setting. For $r>1$, we set
\begin{equation}
\begin{aligned}
V_r &:= \{\vv; \vv\in W^{1,r}(\o; \R^d), \; \diver \vv = \0 \textrm{ in } \o,~\vv \cdot \n =\0 \textrm{ on } \po \},\\
H&:=\overline{V_r}^{L^2(\o; \R^d)} = \{\vec{u} \in E(\Omega);~\diver\vec{u} = 0 \textrm{ in } \Omega,~\vec{u}\cdot\vec{n} = 0 \textrm{ on } \partial\Omega\}, \\ %\left(= L^2_{0,\diver}(\o):=\left\{\vv \in L^2(\Omega;\R^d); \, \int_{\Omega} \vv \cdot \nabla \varphi =0 \textrm{ for all } \varphi\in W^{1,2}(\Omega) \right\}\right),\\
V_r^* &:= (V_r)^*.
\end{aligned}\label{spaces}
\end{equation}
%and equip the space $V_r$ {with} the norm\footnote{ Note that in case $r\in(1, 2d/(d+2))$, we would need to redefine $V_r := \{\vv; \vv\in W^{1,r}(\o; \R^d) \cap L^2(\o;\R^d), \; \diver \vv = \0 \textrm{ in } \o,~\vv \cdot \n =\0 \textrm{ on } \po \}$. However, we work with $r\geq 2d/(d+2)$ and thus we can use that $W^{1,r}(\o) \hookrightarrow  L^2(\o)$.} $\|\vv\|_{V_r}:= \|\nabla \vv\|_{L^r(\o{;\R^{d\times d}})} + \|\vv\|_{L^2(\o{;\R^{d}})}$. For the space $H$, we always write $L^2$-norm (and not $H$-norm).
Here, $E(\Omega):=\{\vec{u}\in L^2(\Omega;\R^d);~\diver\vec{u}\in L^2(\Omega)\}$,  for which it is known that the trace operator has well defined normal component (in the sense of distributions) and there holds $(\vec{u}\cdot \vec{n})\in  (W^{1/2,2}(\partial \Omega; \R^d))^*$, see e.g. \cite{CoFo88}. Referring back to \eqref{spaces}, for any $r\in [\frac{2d}{d+2},\infty)$ and $z>r$, one has
\begin{equation}\label{gelfand}
V_z \hookrightarrow V_r \hookrightarrow H \equiv H^*  \hookrightarrow V_r^* \hookrightarrow V_z^*,
\end{equation}
where all the embeddings are continuous and dense. %Therefore, these spaces form a Gelfand triplet. For simplicity, we also set $V:=V_2$ and $V^*:=V_2^*$. Note that $V$ and $H$ are Hilbert spaces and note that the closure of the space $V_r$ with respect to $L^2$ norm defines the same $H$ for any $r\in [2d/(d+2),\infty)$.  Also, the duality in $V_r$ is defined in a standard way using $\{\f^k\}_{k\in \mathbb{N}}$ a sequence in $H$ converging to $\f$ in $V^*_r$.

Also, we define
\begin{equation}
\begin{aligned}
\cV_r &:= \{\vv; \vv\in W^{1,r}(\o; \R^d), \;~\vv \cdot \n =\0 \textrm{ on } \po \},\\
\cH&:=\overline{\cV_r}^{L^2(\o; \R^d)} = L^2(\o;\R^d),\\
\cV_r^* &:= (\cV_r)^*,
\end{aligned}\label{spaces2}
\end{equation}
and similarly, for $r\in [\frac{2d}{d+2},\infty)$, we have %can define another Gelfand triplet
\begin{equation}\label{gelfand2}
\cV_r \hookrightarrow \cH \equiv \cH^*  \hookrightarrow \cV_r^*.
\end{equation}
Finally, we define
\begin{align*}
\C([0,T];H) &:= \{f \in L^{\infty}(0,T;H); [0,T] \ni t^n \to t \implies f(t^n) \to f(t)~\text{ strongly in } H\}, \\
\C_w([0,T];H) &:= \{f \in L^{\infty}(0,T;H); [0,T] \ni t^n \to t \implies f(t^n) \tow f(t) ~\text{ weakly in } H\}.
\end{align*}

We are now in a position to precisely formulate our main result.
\begin{theorem}\label{genresult}
    Let $\o \in \C^{0,1}$, $T>0$, $\b\in L^{r'}(0,T;V^*_r)$ and $\vv_0 \in H$ be arbitrary. Let $\G$ and $\g$ be arbitrary functions satisfying (G1) - (G4) with $r\in(\frac{2d}{d+2},\infty)$ and (g1)--(g4) with $q\in(1,\infty)$. Then there exists a weak solution to~\eqref{pepa3}, \eqref{pepa9} and \eqref{pepa10} in the following sense: there exist $(\vv, \S, \s)$ such that for $z:= \max\{r, q, \frac{(d+2)r}{(d+2)r -2d}\}$,
    \begin{align*}
        \vv &\in L^r(0,T; V_r) \cap \mathcal{C}_w([0,T];H),\\
        \pt \vv &\in L^{z'}(0,T; V^*_z), \\
        \S &\in L^{r'}(Q), \\
        \s &\in L^{q'}(\Gamma), \\
        \vv &\in L^q(\Gamma),
    \end{align*}
    the balance of linear momentum is satisfied in a weak sense, i.e. for a.a. $t\in(0,T)$ and for all $\vp \in V_z$
    \begin{equation}\label{WF}
        \langle \pt \vv, \vp \rangle_{V_z} - \int_\o (\vov):\nabla \vp \d \x + \int_\o \S:\DD\vp \d \x + \int_{\po} \s \cdot \vp \d S = \langle \b, \vp \rangle_{V_z},
    \end{equation}
    the constitutive equations \eqref{pepa9} and \eqref{pepa10} hold a.e. in $Q$ and $\Gamma$, i.e.,
\begin{align}
\G(\S,\D) &= \0 \qquad \text{ for a.a. } (t,\x) \in Q, \label{pepa9ae} \\
\g(\s,\vvt) &=\0 \qquad \text{ for a.a. } (t,\x) \in \Gamma, \label{pepa10ae}
\end{align}
and the initial condition is attained in the strong sense, i.e.,
    \begin{equation}\label{IC}
        \lim_{t\to0_+} \|\vv(t) - \vv_0\|_2=0.
    \end{equation}
    Also, for all $t\in(0,T)$ the energy inequality holds, i.e.,
    \begin{equation}\label{enineqr}
    \frac{1}{2} \|\vv(t)\|_2^2 + \int_{Q^t} \S:\D \d \x \d \tau + \int_{\Gamma^t} \s\cdot\vv \d S \d \tau \leq \int_0^{t} \langle \b, \vv \rangle_{V_r} \d \tau + \frac12 \|\vv_0\|_2^2.
    \end{equation}
    In addition, if $\o\in \C^{1,1}$ and $\b \in L^{r'}(0,T; \cV_r^*)$, then there exists a pressure $p\in L^{z'}(0,T; L^{z'}(\o))$ such that
   \begin{equation}\label{WFp}
        \langle \pt \vv, \vp \rangle_{\cV_z} - \int_\o (\vov):\nabla \vp \d \x + \int_\o \S:\DD\vp \d \x + \int_{\po} \s \cdot \vp \d S = \langle \b, \vp \rangle_{\cV_z} + \io p \diver \vp \d \x
    \end{equation}
    holds for all $\vp \in \cV_z$ and a.a. $t \in (0,T)$.
\end{theorem}
As stated in Sect.~\ref{Sect2}, the assumption (g4) eliminates no-slip and perfect-slip boundary conditions from the analysis presented here. However, one can incorporate both conditions into the analysis. The case of $\s = \0$ on $\Gamma$ is in fact easy as the boundary term just vanishes. For no-slip boundary conditions, one needs to change the function space for the velocity and pay attention to differences associated with the reconstruction of the pressure (see \cite{wolf} and \cite{BlMaRa2020} for details).

In the proof of Theorem~\ref{genresult}, we use the following powerful convergence result that is a consequence of the properties of the suitably constructed Lipschitz approximations of Bochner functions, see \cite{BDS}. Here, we provide a simplified version (omitting the discussion concerning Lipschitz approximations) suited to the analysis in Sect.~\ref{Sect6}.
\begin{lemma}\label{BDS}
For any interval $I_0\subset (0,T)$ and any ball $B_0\subset\Omega$, set $Q_0:= I_0 \times B_0$. Assume that for $\delta \to 0+$ the following convergences hold:
\begin{align*}
\u^\delta &\tow 0 &&\text{weakly in } L^r(I_0; W^{1,r}(B_0; \R^d)), \\
\u^\delta &\tow^* 0 &&\text{weakly$^*$ in } L^\infty(I_0; L^2(B_0; \R^d)), \\
\u^\delta &\to 0 &&\text{strongly in } L^1(Q_0; \R^d), \\
\H^\delta_1 &\tow 0 &&\text{weakly in } L^{r'}(Q_0; \R^{d\times d}), \\
\H^\delta_2 &\to 0 &&\text{strongly in } L^{1+\e}(Q_0; \R^{d\times d}),
\end{align*}
 whereas $\u^{\delta}$, $\H^{\delta}_1$ and $\H^{\delta}_2$ satisfy
\begin{align*}
\diver\u^{\delta} &= 0 \qquad \textrm{ for a.a. } (t,\x) \in Q_0,  \\
\int_{Q_0} \u^\delta \cdot \pt \vp &- (\H^\delta_1+\H^\delta_2):\nabla \vp \d \x \d t = 0 \textrm{ for all } \vp\in\C^\infty_0(I_0;\C^{\infty}_{0}(B_0; \R^{d})) \textrm{ with } \diver\vp = 0,
\end{align*}
which is a weak formulation of
\begin{equation*}
\pt \u^\delta - \diver (\H^\delta_1+\H^\delta_2) = - \nabla p^\delta.
\end{equation*}
Then there exists a $\xi \in \C^{\infty}_0(Q_0)$ such that
\begin{equation}\label{1816}
\chi_{\frac{1}{8}Q_0} \leq \xi \leq \chi_{\frac{1}{6}Q_0},
\end{equation}
and for every $k \in \N$ there exists a family $\{Q_{\delta,k} \}_{\delta \in (0,1)}$ fulfilling
\begin{equation}\label{qnk}
Q_{\delta,k} \subset Q_0, \quad \limsup_{\delta \to \infty} |Q_{\delta,k}| \leq 2^{-k}
\end{equation}
such that for every $\oS \in L^{r'}(Q_0)$,
\begin{equation}\label{hvi}
\limsup_{\delta \to 0_+} \left| \int_{Q_0} (\H^\delta_1 + \oS)\cdot \nabla \u^\delta \xi \chi_{Q_0 \setminus Q_{\delta,k}} \textrm{d} \x \textrm{d} t\right|\leq C 2^{\frac{-k}{r}}.
\end{equation}
\end{lemma}
\begin{proof}
See \cite[Theorem 2.2 and Corollary 2.4]{BDS}, which is adapted to our setting. Compare also with~\cite{DiRuWo10,BGMS2}.
\end{proof}

\section{Proof of the main result}\label{Sect6} We take any $\o \in \C^{0,1}$, $T>0$, $\b\in L^{r'}(0,T;V^*_r)$ and $\vv_0 \in H$ and fix them for the rest of the proof. Similarly, we consider arbitrary but fixed functions $\G$ and $\g$ satisfying (G1)--(G4) with $r\in(\frac{2d}{d+2},\infty)$ and (g1)--(g4) with $q\in(1,\infty)$. We prove Theorem~\ref{genresult} by means of a two-level approximation.

\subsection{\texorpdfstring{$(\e,\delta)$}{e,d}-approximations} We take two parameters $\e$ and $\delta$ satisfying $\e,\delta \in (0,1)$.

The parameter $\e$ is used to approximate the constitutive equations \eqref{pepa9} and \eqref{pepa10} in the same way as presented in Sect.~\ref{Sect4}, see formulas \eqref{def:Ge} and \eqref{def:ge}. The motivation for such a choice is that the approximations $\Ge$ and $\ge$ possess much better properties in comparison with the properties of $\G$ and $\g$, as summarized in Lemma~\ref{convGg}. In particular, for $\e$-approximation,  we can use \eqref{bddSD} and thus stay on the level of nonlinear yet uniformly Lipschitz continuous and uniformly monotone operators. Consequently, we come from the $(r,r')$-coercivity for $\D$ and $\S$ to the $L^2$-coercivity for the $\e$-approximations $\De$ and $\Se$. The mathematical theory for problems of parabolic type with such nonlinear operators is well known, see \cite[Appendix~C]{BMM21} for example. Since this type of $\e$-approximation changes the $r$-structure to $2$-structure, it is suitable to define the following auxiliary numbers
\begin{align*}
&\mu := \min\{r,2\},  &&\mu' := \max\{r',2\}, \\
&\nu := \min\{r',2\}, &&\nu' := \max\{r,2\}.
\end{align*}

The purpose of introducing another approximation parameter $\delta$ is twofold. First,
we need to regularize the right-hand side $\b$. Using a density argument, we approximate $\b \in L^{r'}(0,T;V_r^*)$ by a sequence of $\b^\delta \in L^{\mu'}(0,T; V_{\mu}^*)$ such that
\begin{equation}\label{bstro}
\b^\delta \to \b \qquad \textrm{ strongly in } L^{r'}(0,T;V_r^*) \textrm{ as } \delta \to 0_+.
\end{equation}
Second, $\delta$ will be used in the definition of the cut-off function $\phi_{\delta}$ used to control (bound) the convective term and the influence of the boundary condition. More precisely,
%\subsection{\texorpdfstring{$(\e,\delta)$}{e,d} - approximating problem}
%For arbitrary $\e,\delta \in (0,1)$ we start with the following $(\e,\delta)$-approximating problem:
%\begin{equation}
%\label{2}
%\begin{aligned}
%&\text{For given  } \o ~\text{Lipschitz},~T>0,~\vv_0: \o \to \R^d,~\b: Q \to \R^d,
%\\ \Ge(\S,\D) &:= \G(\S-\e\D, \D-\e\S) \text{ for } \G \text{ satisfying }\text{(G1)--(G4) with } r>\frac{2d}{d+2},
%\\ \ge(\s,\vvt) &:= \g(\s-\e\vvt, \vvt-\e\s) \text{ for } \g \text{ satisfying }\text{(g1)--(g4) with } q>1,\\ &\qquad\text{find a triplet } (\vv, \S, \s) \text{ a weak solution to}\\
%&\begin{aligned}
%\diver \vv =&~ 0 &&\text{in } Q, \\
%\pt \vv+ \diver [(\vov)\phi_{\delta}(|\vv|^2)] =&~ -\nabla p + \diver \S + \b &&\text{in } Q, \\
%\Ge(\S,\D) =&~\0 &&\text{in } Q, \\
%\ge(\s, \vv_\tau) =&~ \0 &&\text{on } \Gamma, \\
%\vv \cdot \n =&~0 &&\text{on } \Gamma, \\
% \vv (0,\cdot) =&~\vv_0 &&\text{in } \o,
%\end{aligned}\\
%\text{ where } &-(\S \n )_{\tau} =: \phi_{\delta}(|\vv|^2) \s \text{ and } \phi_\delta \text{ is defined in}~\eqref{phid}.
%\end{aligned}
%\end{equation}
$\phi_{\delta}$ is defined through
\begin{equation}\label{phid}
    \phi_\delta(s):=\phi(\delta s)  \hspace{1cm} \text{ where } \hspace{1cm}
    \phi(s) := \begin{cases}
        1   &|s|\in[0,1),\\
        2-s &|s|\in[1,2),\\
        0   &|s|\in[2,\infty).
    \end{cases}
\end{equation}
Note that, for all $s\geq 0$, $\phi_\delta(s) \to 1$ as $\delta \to 0_+$.

Thus, the role of such a $\delta$-approximation is to avoid possible singularities in lower order terms. This together with the $\e$-approximation turns the problem into one in the classical setting of monotone operators with compact lower order terms. Consequently, we state without proof the following theorem concerning the existence of a weak solution to the $(\e,\delta)$-approximating problem.
\begin{lemma}\label{app}
 Let $\e,\delta\in (0,1)$ be arbitrary but fixed. There exist $(\vv^{\e,\delta}, \S^{\e,\delta},\s^{\e,\delta})$ solving the $(\e,\delta)$-approximation of the problem \eqref{pepa3}, \eqref{pepa9} and \eqref{pepa10} in the following sense:
    \begin{align*}
        \vv^{\e,\delta} &\in L^2(0,T; V_2) \cap \mathcal{C}([0,T];H),\\
        \pt \vv^{\e,\delta} &\in L^{2}(0,T; V_2^*), \\
        \S^{\e,\delta} &\in L^{2}(Q), \\
        \s^{\e,\delta} &\in L^2(\Gamma);
    \end{align*}
    the balance of linear momentum is satisfied in a weak sense, i.e.
    \begin{equation}\label{WFed}
    \begin{aligned}
        \langle \pt \vv^{\e,\delta}, \vp \rangle_{V_2} &- \!\int_\o \phi_{\delta}(|\vv^{\e,\delta}|^2) (\vv^{\e,\delta} \otimes \vv^{\e,\delta})\!:\!\nabla \vp \d \x + \int_\o \S^{\e,\delta}\!:\!\DD\vp \d \x \\
        &+ \!\int_{\po} \phi_{\delta}(|\vv^{\e,\delta}|^2) \s^{\e,\delta} \!\cdot \vp \d S = \langle \b^\delta, \vp \rangle_{V_2} ~\text{ for a.a. } t\in(0,T) \text{ and all } \vp \in V_2;
    \end{aligned}
    \end{equation}
    the constitutive equations \eqref{def:Ge} and \eqref{def:ge} are fulfilled pointwise, i.e.,
    \begin{align}
    \Ge(\S^{\e,\delta},\DD\vv^{\e,\delta}) &= \0 \quad \textrm{ for a.a. } (t,\x)\in Q \label{pepa18}\\ \ge(\s^{\e,\delta},\vv^{\e,\delta}_\tau) &=\0 \quad \textrm{ for a.a. } (t,\x) \in \Gamma, \label{pepa19}
    \end{align}
    and $\vv^{\e,\delta}(0,\cdot) = \vv_0$ in $\o$; as $\vv^{\e,\delta}\in\mathcal{C}([0,T];H)$, this means that $\vv_0$ is attained strongly.
\end{lemma}
%We already described that Theorem~\ref{app} can be proven by using techniques which are standard by now and will not be discussed here. On the other hand, we focus on the passage of the limit $\e\to 0_+$ and $\delta \to 0_+$ to prove the existence of a weak solution to the general problem~\eqref{pepa3}, \eqref{pepa9} and \eqref{pepa10}.

\subsection{Limit \texorpdfstring{$\e \to 0_+$}{to}}\label{subsube}
{ In this part, we keep $\delta\in (0,1)$ fixed and we let $\e\to 0_+$. We also write $(\ve,\Se,\se)$ instead of $(\vv^{\e,\delta}, \S^{\e,\delta},\s^{\e,\delta})$.} %Since this limit changes the $2$-structure to $r$-structure, we define the auxiliary numbers
%\begin{align*}
%&\mu := \min\{r,2\},  &&\mu' := \max\{r',2\}, \\
%&\nu := \min\{r',2\}, &&\nu' := \max\{r,2\}.
%\end{align*}
%Note that $\b^\delta \in L^{\mu'}(0,T; V_{\mu}^*)$. Since $\mu \le 2$, we can use Lemma~\ref{app} and construct for every $(\e,\delta)$ a solution to \eqref{WFed}.}

\subsubsection{Uniform \texorpdfstring{$\e$}{e}-independent estimates}\label{subcislo}
To %proceed with the limit $\e\to 0_+$, we first
obtain a~priori estimates which are uniform with respect to $\e$ (but may depend on $\delta$), we set $\vp := \ve$ in~\eqref{WFed}, integrate the result over $(0,t)$, and use the facts that $\pt \ve \in L^2(0,T;V_2^*)$ and $\int_\o \fidve (\vove) : \nabla \ve \d \x = 0$\footnote{Let $\Phi_\delta$ denote the primitive function to $\phi_\delta$. Using integration by parts we observe that
$$
2\io \fidve (\vove) : \nabla \ve \d \x =  \io \fidve \ve \cdot \nabla|\ve|^2\d \x =   \io \ve \cdot \nabla \Phi_\delta (|\ve|^2)\d \x = \io \diver \left(\Phi_\delta (|\ve|^2) \ve  \right) \d \x =0,
$$
referring to Gauss' theorem and the boundary condition $\ve\cdot\n = 0$ on $\partial \Omega$ in the last identity.}. This gives
\begin{equation*}
    \frac12 \|\ve(t)\|_2^2 + \int_{Q^t} \Se:\De \d \x \d \tau + \int_{\Gamma^t} \fidve \se \cdot \ve_\tau \d S \d \tau = \int_0^t \langle \b^\delta,\ve\rangle_{V_\mu} \d \tau + \frac12 \|\vv_0\|_2^2.
\end{equation*}
As $(\Se,\De)$ are the null points of $\Ge$ and $(\se,\ve_\tau)$ are the null points of $\ge$, by applying ~\eqref{indeG} and~\eqref{indeg} from  Lemma~\ref{convGg} with one simple estimation applied to~\eqref{indeg}, we obtain
\begin{align}
    \Se:\De &\geq \tilde{C}_1(|\Se|^\nu + |\De|^\mu) - \tilde{C}_2, \label{pepa22}\\
    \se\cdot \ve_\tau &\geq \hat{c}_1 (|\se| +|\ve_\tau| )^{\min\{q, q'\}} - \hat{c}_2. \label{pepa23}
\end{align}
Hence
\begin{equation}\label{wfve}
\begin{aligned}
\frac12 \|\ve(t)\|_2^2 &+ \tilde{C}_1 \iqt |\Se|^{\nu} + |\De|^{\mu} \d \x \d \tau + \hat{c}_1 \igt \fidve(|\se| + |\ve_\tau| )^{\min\{q, q'\}}\d S \d \tau  \\
&\leq \frac12 \|\ve(t)\|^2_2 + \iqt \Se : \De \d \x \d \tau + \igt \fidve \se \cdot \ve_\tau \d S \d \tau +C\\
&\leq \it \langle \b^\delta,  \ve \rangle_{V_\mu} \d \tau + \frac12 \|\vv_0\|_2^2 +C.
\end{aligned}
\end{equation}
To estimate the term with $\b^{\delta}$, we incorporate Korn's and Young's inequalities and conclude that
\begin{equation}\label{pepa24}
\begin{aligned}
\langle \b^\delta,  \ve \rangle_{V_\mu} &\leq \|\b^\delta\|_{V_\mu^*} \left(\|\ve\|_2 + \|\De \|_{\mu}\right) \\
&\leq \frac{\tilde{C}_1}{2} \|\De \|_{\mu}^\mu + C\left(\|\b^\delta\|_{V_\mu^*}^{\mu'} + (1+ \|\ve\|_2^2) \|\b^\delta\|_{V_\mu^*} \right).
\end{aligned}
\end{equation}
Inserting this estimate into \eqref{wfve} and applying Gronwall's lemma we obtain
\begin{equation}\label{venekH}
\sup_{t \in (0,T)} \|\ve(t)\|_2 \leq C(\delta).
\end{equation}
With this information we look at \eqref{wfve} and \eqref{pepa24} again and
conclude that
%the see that the last term Together, taking supremum over $t \in (0,T)$ in %\eqref{wfve} and using assumption on data ($\b^\delta \in L^{\mu'}(0,T;V_{\mu}^*)$ and $\vv_0 \in H$) we get
\begin{equation}\label{unif}
\begin{aligned}
\sup_{t \in (0,T)} \|\ve(t)\|^2_2 &+ \iq |\Se|^{\nu} + |\De|^{\mu} \d \x \d \tau\leq C(\delta),
\end{aligned}
\end{equation}
which gives
\begin{equation}\label{vemumu1}
\|\ve\|_{L^\mu(0,T;V_{\mu})} \leq C(\delta).
\end{equation}
Also, again from \eqref{wfve} and \eqref{pepa24},
\begin{equation}\label{SeDveC}
\iq \Se : \De \d \x \d \tau \leq C(\delta).
\end{equation}
Finally, applying Young's inequality to the left-hand side of \eqref{pepa23}
we get
$$
\hat{c}_1 (|\se| + |\ve_\tau|)^{\min\{q, q'\}} \le \frac{\hat{c}_1}{2} |\se|^{\min\{q, q'\}} + c_3 |\ve_\tau|^{\max\{q, q'\}} + \hat{c}_2,
$$
which gives
\begin{equation}\label{pepa25}
\frac{\hat{c}_1}{2} |\se|^{\min\{q, q'\}} \le c_3 |\ve_\tau|^{\max\{q, q'\}} + \hat{c}_2.
\end{equation}
Multiplying this inequality by $[\fidve]^{\min\{q, q'\}}$ and using then the truncation properties of $\fidve$ introduced in~\eqref{phid} we get
\begin{equation}\label{omezse}
\left[\fidve |\se|\right]^{\min\{q, q'\}} \leq C [\fidve]^{\min\{q, q'\}} (1+|\ve_\tau|^{\max\{q, q'\}}) \leq C(\delta)
\end{equation}
Hence, \eqref{omezse} imply that
\begin{equation}\label{vemumu2}
\|\fidve \se\|_{L^\infty (\Gamma)} \leq C(\delta).
\end{equation}

To estimate the time derivative of $\ve$, we set $\du := \{ \w \in V_r \cap V_2; \|\w\|_{V_{\nu'}} \leq 1\}$. Since $\du \subset V_2$ we can set $\vp := \w \in \du$ in the equation \eqref{WFed} and conclude that (using~\eqref{omezse})
\begin{align*}
\|\pt \ve\|_{V_{\nu'}^*} &= \sup_{\du} \langle \pt \ve, \w \rangle_{V_{\nu'}} \\
&= \sup_{\du} \left(\io \fidve (\vove):\nabla \w \d \x -\io \Se : \DD \w \d \x - \right. \\
&\qquad \hspace{3cm}\left. \ipo \fidve( \se \cdot \w )\d S + \langle \b^\delta,  \w \rangle_{V_{\mu}} \right) \\
&\leq \sup_{\du} \left( C(\delta) \|\w\|_{\du} + \|\Se\|_{\nu} \|\DD \w\|_{\nu'} + \|\b^\delta\|_{V_{\mu}^*} \|\w\|_{V_{\mu}} \right).
\end{align*}
As $V_{\mu} \subset V_{\nu'}$, raising the last inequality to the power $\nu$ and integrating the result over $(0,T)$ we obtain, using also \eqref{unif},
\begin{equation}\label{unifpt}
\begin{aligned}
\iT \|\pt \ve\|^\nu_{V_{\nu'}^*} \d t &\leq \iT \|\Se\|^\nu_{\nu} + \|\b^\delta\|^\nu_{V_{\mu}^*} \d t \leq C(\delta).
\end{aligned}
\end{equation}

\subsubsection{Limit passage \texorpdfstring{$\e\to 0_+$}{e}}

The estimates \eqref{venekH}, \eqref{unif}, \eqref{vemumu1}, \eqref{vemumu2} and \eqref{unifpt} imply the existence of subsequences (that we label again as the original sequences) and the limiting objects such that as $\e \to 0_+$
\begin{equation}\label{converg}
\begin{aligned}
\ve &\tow \vv &&\text{weakly in } L^{\mu}(0,T;V_{\mu}), \\
\Se &\tow \S &&\text{weakly in } L^{\nu}(Q), \\
\ve &\tow^* \vv &&\text{weakly$^*$ in } L^{\infty}(0,T;H), \\
\fidve\se &\tow^* \overline{\fidv\s} &&\text{weakly$^*$ in } L^\infty(\Gamma), \\
\pt \ve &\tow \pt \vv &&\text{weakly in } L^{\nu}(0,T;V_{\nu'}^*).
\end{aligned}
\end{equation}
By the Aubin--Lions Compactness Lemma and the Trace Theorem we also observe
\begin{equation}\label{strongv}
\begin{aligned}
  \ve &\to \vv &&\text{strongly in } L^{2}(Q), \\
  \ve &\to \vv &&\text{strongly in } L^{1}(\Gamma).
\end{aligned}
\end{equation}
Consequently,
\begin{equation}\label{strongvov}
    \fidve (\vove) \to \fidv (\vov) ~~~\text{ strongly in } L^{s}(Q)~\text{ for any } s\in [1, \infty).
\end{equation}
Furthermore, from~\eqref{strongv}, for any $\sigma>0$ there exists $\Gamma_\sigma \subset \Gamma$ such that $|\Gamma \setminus \Gamma_\sigma|\leq \sigma$ and
\begin{equation}\label{strongvgammasigma}
   \ve \to \vv ~~~\text{ strongly in } L^{\infty}(\Gamma_\sigma) \qquad\textrm{ (modulo subsequence)}
\end{equation}
As a consequence of~\eqref{strongvgammasigma} and \eqref{pepa25}, we get
\begin{equation}\label{weaksgammasigma}
   \se \tow^* \s ~~~\text{ weakly$^*$ in } L^{\infty}(\Gamma_\sigma).
\end{equation}
Moreover, \eqref{SeDveC} in combination with  Lemma~\ref{convGg} implies that
\begin{equation}\label{SerDer}
\S \in L^{r'}(Q) ~\text{ and }~ \D \in L^r(Q).
\end{equation}

Now, we integrate~\eqref{WFed} over $(0,T)$ and study the limit $\e \to 0_+$ using  the convergence results~\eqref{converg}--\eqref{weaksgammasigma}. For any $\w \in L^{\nu'}(0,T;V_{\nu'})$ we end up with
\begin{equation}\label{WFt}
\begin{split}
&\iT \langle \pt \vv,  \w \rangle_{V_\mu} \d \tau + \iq \S : \DD \w \d \x \d \tau + \ig \overline{\fidv \s} \cdot \w \d S \d \tau \\
&\qquad = \iT \langle \b^\delta, \w \rangle_{V_\mu} \d \tau +\iq \fidv (\vov):\nabla \w \d \x \d \tau %\qquad \textrm {for any } \w \in L^{\nu'}(0,T;V_{\nu'}).
\end{split}
\end{equation}
Thanks to the dense embedding $V_r \cap V_2 \hookrightarrow V_r$, we directly obtain
\begin{equation*}
\iT \|\pt \vv\|^{r'}_{V_r^*} \d t \leq \iT \|\S\|^{r'}_{r'} + \|\fidv\vov\|^{r'}_{r'}  + \|\overline{\fidv \s} \|^{r'}_{L^{r'}(\po)} + \|\b^\delta\|^{r'}_{V_r^*}\d t \leq C(\delta),
\end{equation*}
where the last inequality follows from the properties of $\fidv$,~\eqref{SerDer} and~\eqref{vemumu2}.

Moreover, thanks to $\vv \in L^r(0,T;V_r)$, $\pt \vv \in L^{r'}(0,T;V_r^*)$, and the Gelfand triple \eqref{gelfand}, there holds $\vv \in \C([0,T];H)$. In addition, in a standard way we can show that $\vv(0,\cdot) = \vv_0$.

\subsubsection{Identification of nonlinearities}

To complete the limit $\e \to 0+$, we need to verify that $\G(\S(t,\x),\D(t,\x))=\0$ for a.a. $(t,\x)\in Q$, $\g(\s(t,\x),\vvt(t,\x))=\0$ for a.a. $(t,\x)\in \Gamma$ and  $\overline{\fidv \s}= \fidv \s$ on $\Gamma$.

We first observe that for all $t \in (0,T)$,
\begin{equation}\label{goal1}
\limsup_{\e \to 0_+} \iqt \Se : \De \d \x \d \tau + \igt \fidve \se \cdot \ve \d S \d \tau \leq \iqt \S : \D \d \x \d \tau + \igt \overline{\fidv \s} \cdot \vv \d S \d \tau.
\end{equation}
Indeed, setting $\vp:=\ve$ in \eqref{WFed} and integrating the result over $(0,t)$ for $t\in(0,T)$, we get
\begin{equation*}
\iqt \Se : \De \d \x \d \tau + \igt \fidve \se \cdot \ve \d S \d \tau = \it \langle \b^\delta,  \ve \rangle_{V_\mu} \d \tau + \frac12 \|\vv_0\|_2^2- \frac12 \|\ve(t)\|^2_2.
\end{equation*}
Taking the $\limsup_{\e\to 0+}$ and using the  weak lower semicontinuity of the $L^2$-norm w.r.t. weakly* converging sequence in $L^{\infty}(0,T;H)$, we obtain for a.a. time
\begin{equation}\label{slabave}
\begin{aligned}
\limsup_{\e \to 0_+} &\iqt \Se : \De \d \x \d \tau + \igt \fidve\se \cdot \ve \d S \d \tau \\
&\leq \it \langle \b^\delta,  \vv \rangle_{V_\mu} \d \tau + \frac12 \|\vv_0\|_2^2- \frac12 \|\vv(t)\|^2_2.
\end{aligned}
\end{equation}
On the other hand, we can set $\w:=\vv \chi_{(0,t)}$ in \eqref{WFt} (we already have the right duality pairings to do so) and integrate the result over $(0,t)$. We conclude that
\begin{equation}\label{slabav}
\iqt \S : \D \d \x \d \tau + \igt \overline{\fidv \s} \cdot \vv \d S \d \tau = \it \langle \b^\delta,\vv \rangle_{V_\mu} \d \tau + \frac12 \|\vv_0\|_2^2- \frac12 \|\vv(t)\|^2_2.
\end{equation}
Comparing \eqref{slabav} and \eqref{slabave} we see that~\eqref{goal1} holds true.

Finally, using \eqref{goal1} we want to prove that for a.a. $t \in (0,T)$,
\begin{align}
\limsup_{\e \to 0_+} \iqt \Se : \De \d \x \d \tau &\leq \iqt \S : \D \d \x \d \tau \label{SD}\\
\limsup_{\e \to 0_+} \int_{\Gamma_\sigma} \se \cdot \ve \d S \d \tau &\leq \int_{\Gamma_\sigma} \s \cdot \vv \d S \d \tau.\label{sv}
\end{align}
where $\Gamma_\sigma$ is introduced around ~\eqref{strongvgammasigma}. However,~\eqref{sv} directly follows from~\eqref{strongvgammasigma} and~\eqref{weaksgammasigma} (in fact, with the equality sign in \eqref{sv}). Recalling  ~\eqref{bddseve} in Lemma~\ref{convGg}, \eqref{sv} implies that  $\g(\s,\vvt)=\0$ a.e. in $\Gamma_\sigma$. However, since $\sigma>0$ can be made arbitrarily small, the statement $\g(\s,\vvt)=\0$ holds true  a.e. on $\Gamma$ (again considering a suitably chosen subsequence).

Next, we show that
\begin{equation}\label{odstrpruh}
    \liminf_{\e\to 0_+} \igt \fidve \se \cdot \ve \d S \d \tau  \geq  \igt \overline{\fidv \s} \cdot \vv \d S \d \tau.
\end{equation}
First, from~\eqref{strongvgammasigma} and~\eqref{weaksgammasigma} it follows that $ \overline{\fidv \s}= \fidv \s$ on $\Gamma_\sigma$. Then thanks to~\eqref{nonnegge} it follows that
\begin{align*}
     \liminf_{\e\to 0_+} &\igt \fidve \se \cdot \ve \d S \d \tau \geq  \liminf_{\e\to 0_+} \int_{\Gamma^t \cap \Gamma_\sigma} \fidve \se \cdot \ve \d S \d \tau \\
     &= \int_{\Gamma^t \cap \Gamma_\sigma} \fidv \s \cdot \vv \d S \d \tau \\
     &= \int_{\Gamma^t \cap \Gamma_\sigma} \overline{\fidv \s} \cdot \vv \d S \d \tau .
\end{align*}
Letting $\sigma \to 0_+$ in the above inequality we obtain~\eqref{odstrpruh}.

Using~\eqref{odstrpruh} in~\eqref{goal1} we obtain~\eqref{SD}, which is the assumption~\eqref{bddSeDe} in Lemma~\ref{convGg}. Consequently, $\G(\S,\D)=\0$ a.e. in $Q^t$ for a.a. $t$.

We summarize the results proved above (in Subsection~\ref{subsube}) in the following lemma.
\begin{lemma}\label{appd}
    % Let $\delta\in (0,1)$ be arbitrary and $\o \subset \R^d$ be a Lipschitz domain, $T>0$, $\b^\delta\in L^{\mu'}(0,T;V_{\mu}^*)$ and $\vv_0 \in H$. Let $\G$ satisfy (G1)--(G4) for $r > \frac{2d}{d+2}$ and $\g$ satisfy (g1)--(g4) for $q>1$. Then
For any $\delta \in (0,1)$ there exists a triplet $(\vv^\delta, \S^\delta, \s^\delta)$ such that
    \begin{align*}
        \vv^\delta &\in L^r(0,T; V_r) \cap \mathcal{C}([0,T];H),\\
        \pt \vv^\delta &\in L^{r'}(0,T; V_r^*), \\
        \S^\delta &\in L^{r'}(Q), \\
        \s^\delta &\in L^{q'}(\Gamma);
    \end{align*}
    \begin{equation}\label{WFd}
        \begin{split}\langle \pt \vv^\delta, \vp \rangle_{V_r} &- \int_\o \phi_{\delta}(|\vv^\delta|^2)(\vv^\delta \otimes \vv^\delta):\nabla \vp \d \x + \int_\o \S^\delta:\DD\vp \d \x  + \int_{\po} \phi_{\delta}(|\vv^\delta|^2) \s^\delta \cdot \vp \d S \\& = \langle \b^\delta, \vp \rangle_{V_r} \quad \textrm{ for a.a. } t\in(0,T) \textrm{ and for all } \vp \in V_r,\end{split}
    \end{equation}
    \begin{align}
        \G(\S^\delta,\D^\delta) &= \0 \quad \textrm{ a.e. in } Q,\label{pepa30}\\
        \g(\s^\delta,\vv^\delta_\tau)&=\0 \quad \textrm{ a.e. on } \Gamma, \label{pepa31}
    \end{align} and the initial condition $\vv_0$ is attained in the strong sense.
%\begin{equation}\label{icd}
%        \lim_{t\to 0_+} \|\vv(t) - \vv_0\|_2=0.
%    \end{equation}
\end{lemma}

\subsection{Limit \texorpdfstring{$\delta\to 0_+$}{d}}\label{subposl}
We recall that $\b^\delta$ satisfies \eqref{bstro}.
\subsubsection{Uniform (\texorpdfstring{$\delta$}{d}-independent) estimates and their consequences} Setting $\vp := \vd$ in \eqref{WFd} (the convective term vanishes, see  Sect.~\ref{subcislo} for details) we get for a.a. $t\in (0,T)$
\begin{equation}\label{priprava}
\dt \|\vd(t)\|_2^2 + 2\io \Sd \!:\! \Dd \d \x  + 2\ipo \fidvd \sd \!\cdot\! \vd \d S =2\langle \b^{\delta}, \vd\rangle_{V_r}
\end{equation}
Integrating \eqref{priprava} over $(0,t)$, using (G4) and (g4), the assumptions on $\vv_0$ and on the right-hand sides $\b^\delta$ and $\b$ (see \eqref{bstro}), H\"older's, Young's and Gronwall's inequalities, we obtain
\begin{equation*}
\sup_{t \in (0,T)}\|\vd(t)\|_2^2 + \iq \Sd \!:\! \Dd \d \x \d \tau + \ig \fidvd \sd \!\cdot\! \vd \d S \d \tau\leq C.
\end{equation*}
Using (G4) and (g4) and proceeding as in Sect.~\ref{subcislo},  we conclude that
\begin{subequations}\label{uniesti}
    \begin{align}
        \|\vd\|_{L^\infty(0,T;H)\cap L^r(0;T;V_r)} &\leq C, \label{unidvd}\\
\|\Sd\|_{L^{r'}(Q)} &\leq C, \label{uniSd}\\
 \ig \fidvd (|\sd|^{q'} + |\vd|^q) \d S \d \tau &\leq C.\label{unisdvd}
    \end{align}
Thanks to the fact that $0\leq \fidvd\leq 1$, it directly follows from the above estimates that
\begin{equation}\label{unifidvdsd}
    \ig |\fidvd \sd|^{q'} \d S \leq \ig \fidvd|\sd|^{q'} \d S \leq  C.
\end{equation}

Next, we explain the definition of $z$ in Theorem~\ref{genresult}, which is related to the uniform estimate for $\partial_t \vd$. First, it follows from the definition that $z':= \min\{r', q', \frac{(d+2)r}{2d}\}$ and each number in the bracket is the dual exponent to the integrability exponent of the terms on the left-hand side of \eqref{WFd}. We start with the convective term. Thanks to the assumption $r > \frac{2d}{d+2}$, we have  $\frac{(d+2)r}{2d}>1$ and we obtain that
\begin{equation}\label{hbasu}
\iT \|(\vd \otimes \vd) \phi_\delta(|\vd|^2)\|^{\frac{(d+2)r}{2d}}_{\frac{(d+2)r}{2d}} \d t \leq \iT \|\vd\|^{\frac{(d+2)r}{d}}_{\frac{(d+2)r}{d}} \d t \leq \iT \|\vd\|^{\frac{2r}{d}}_2 \|\vd\|^r_{V_r} \d t \leq  C,
\end{equation}
where for the first inequality we used the fact that $\fidvd \leq 1$, for the second inequality we used the standard interpolation in Lebesgue and Sobolev spaces and for the last inequality we used~\eqref{unidvd}.
Then using~\eqref{hbasu},~\eqref{uniSd},~\eqref{unifidvdsd} and the assumption on $\b^{\delta}$ we deduce from~\eqref{WFd} the following estimate
\begin{equation}\label{uniptvd}
\iT \|\pt\vd\|^{z'}_{V_z^*} \d t \leq C \iT \left( \|\Sd\|_{r'} + \|\vd\|^2_{\frac{(d+2)r}{d}} +\|\fidvd\sd\|_{L^{q'}(\po)} + \|\b^{\delta}\|_{V_r^*}\right)^{z'} \d t\leq  C.
\end{equation}
\end{subequations}
Finally, from the uniform estimates~\eqref{uniesti}, the Aubin--Lions lemma and the Trace Theorem we conclude that there are  subsequences of $\{\vd,\Sd,\sd\}$ denoted again by $\{\vd,\Sd,\sd\}$ such that
\begin{subequations}\label{konve}
\begin{align}
\vd &\tow^* \vv &&\text{weakly$^*$ in } L^\infty(0,T;H), \label{vvlih}\\
\vd &\tow \vv &&\text{weakly in } L^r(0,T;V_r),\label{vvvr} \\
(\vd \otimes \vd) \phi_\delta(|\vd|^2) &\to \vv \otimes \vv &&\text{strongly in }~ L^{\rho}(Q) ~\text{ for }~\rho \in \left[1, \frac{(d+2)r}{2d}\right), \label{vdovd}\\
\pt \vd &\tow \pt \vv &&\text{weakly in } L^{z'}(0,T;V_z^*), \label{ptvniekde} \\
\vd &\to \vv &&\text{strongly in } L^r(0,T;L^2(\o)), \\
\vd &\to \vv &&\text{strongly in } L^\gamma(Q) \text{ for } \gamma\in \left[1,\frac{(d+2)r}{d}\right), \label{strong53}\\
\Sd &\tow \S &&\text{weakly in } L^{r'}(Q), \label{Sconv}\\
\fidvd \sd &\tow \s &&\text{weakly in } L^{q'}(\Gamma), \label{shran}\\
\fidvd\vd &\tow \vv &&\text{weakly in } L^{q}(\Gamma), \\
\vd &\to \vv &&\text{strongly in } L^1(\Gamma). \label{silnehran}
\end{align}
\end{subequations}
Then, we consider $\vp \in L^z(0,T;V_z)$ in \eqref{WFd}, integrate the result over $(0,T)$, and apply the convergence results from \eqref{konve}. We obtain
\begin{equation}
\iT \langle \pt \vv, \vp \rangle_{V_z} \d t + \iq (\S-(\vv \otimes \vv)): \nabla \vp \d \x \d t + \ig \s \cdot \vp \d S \d t = \iT \langle \b,  \vp \rangle_{V_z} \d t. \label{zatnizubyatahni}
\end{equation}
Therefore, the weak formulation \eqref{WF} holds for a.a. $t \in (0,T)$. Moreover, the results \eqref{vvlih}, \eqref{vvvr}, and \eqref{ptvniekde} imply that $\vv \in \C_w([0,T];H)$.

Next, we need show that $\G(\S,\D)=0$ a.e. in $Q$ and $\g(\s,\vvt)=0$ a.e. on $\Gamma$. Here, we have several possibilities. First, if $z=r$, i.e. if $q\leq r$ and $r\geq \frac{3d+2}{d+2}$, one can simply  use $\vp:=\vv$ in \eqref{zatnizubyatahni} and therefore, one might mimic the theory developed in~\cite{BMM21}. Next, if $r\geq \frac{3d+2}{d+2}$ but $q>r$, one may observe that the choice $\vp:=\vv$ is admissible in the second and third integral on the left-hand side and also in the term on the right-hand side. Consequently, one might try to generalize the concept of the Gelfand triple and define properly a duality pairing $\langle \pt \vv, \vp\rangle$ and again to mimic the theory from \cite{BMM21}. Finally, in the case $r<\frac{3d+2}{d+2}$ we cannot set $\vp:=\vv$ in the second term and therefore the theory from~\cite{BMM21} cannot be adapted directly to our case. As the novel method developed here, which is based also on the use of Lemma~\ref{BDS}, covers also the simple cases discussed above, we just present a unified  procedure for all values of $r>\frac{2d}{d+2}$.

\subsubsection{Identification of nonlinearities  on the boundary}

By virtue of \eqref{silnehran} and Egoroff's theorem, for every $\eta>0$ there exists $\Gamma_{\eta}$ so that $|\Gamma \setminus \Gamma_{\eta}|<\eta$ and
\begin{equation}\label{pepa100a}
\vd \to \vv  \quad \textrm{ strongly in } L^\infty (\Gamma_{\eta}).
\end{equation}
It follows from (g4) and Young's inequality that
\begin{equation*}
  c_1|\sd|^{q'} \le \sd\cdot\vd + c_2 \leq \frac{c_1}{2} |\sd|^{q'} + c_3 |\vd|^q + c_2 \quad\implies\quad  |\sd|^{q'} \leq C (1+ |\vd|^q).
\end{equation*}
This together with \eqref{pepa100a} implies that $\{\sd\}$ is a bounded sequence on $\Gamma_\eta$.
%and consequently,
%\begin{equation*}
%    \sd \tow^* \s \quad~\text{ weakly$^*$ in } L^{\infty}(\Gamma_{\eta}).
%\end{equation*}
By \eqref{shran} and Lebesgue's Dominated Convergence Theorem, we conclude that, as $\delta \to 0_+$,
\begin{equation*}
\int_{\Gamma_{\eta}} \sd \cdot \vd \d \x \d t  = \int_{\Gamma_{\eta}} \fidvd \sd \cdot \vd \d \x \d t + \int_{\Gamma_{\eta}} (1- \fidvd) \sd \cdot \vd \d \x \d t \to \int_{\Gamma_{\eta}} \s \cdot \vv \d \x \d t.
\end{equation*}
Then from Lemma~\ref{convGg}, $\g(\s,\vv)=\0$ a.e. on $\Gamma_{\eta}$, and letting $\eta \to 0_+$, we obtain that $\g(\s,\vv)=\0$ a.e. on $\Gamma$, and also that for all $\eta$,
\begin{equation}\label{slabasv1}
\sd \cdot \vd \tow \s \cdot \vv \text{ weakly in } L^1(\Gamma_{\eta}).
\end{equation}

\subsubsection{Identification of nonlinearities inside the domain}

Identification in $Q$ is not so straightforward, especially due to the lack of proper duality pairing in the convective term and consequently potential failure of the energy equality for the limiting equation. We start with subtracting the weak formulation for $\vd$ \eqref{WFd} from the one for $\vv$ \eqref{WF} and integrate the difference over $(0,T)$. We deduce that, for all  $\vp \in L^z(0,T;V_z)$,
\begin{equation}\label{dnes}
\begin{aligned}
&\iT \langle \pt (\vd -\vv), \vp \rangle_{V_z} \d t - \iq \left((\vd \otimes \vd) \phi_\delta(|\vd|^2)- \vv \otimes \vv\right): \nabla \vp \d \x \d t  \\
&\qquad + \iq (\Sd -\S): \DD \vp \d \x \d t + \ig (\fidvd\sd- \s) \cdot \vp \d S \d t - \int_0^T \langle \b^{\delta}-\b, \vp\rangle_{V_r} \d t = 0.
\end{aligned}
\end{equation}
Next, we localize the above formulation and also omit writing the duality pairing in $V_r$. Indeed, by using the classical theory for $r$-Stokes problems, we can find\footnote{We can set $\F$ as $\F:=|\nabla \w|^{r-2}\nabla \w$, where $\w$ solves the homogeneous Dirichlet problem
$$
-\diver (|\nabla \w|^{r-2}\nabla \w) = -\nabla \pi+\b, \quad \diver \w=0 \qquad  \textrm{ in } \Omega.
$$
Similarly from $\b^{\delta}$ we come to $\F^{\delta}$.}
$\F^{\delta}$ and $\F$ such that
\begin{equation}\label{FD}
\F^{\delta} \to \F \textrm{ strongly in }L^{r'}(Q)
\end{equation}
that fulfils
$$
\int_0^T \langle \b^{\delta}-\b, \vp\rangle_{V_r} \d t=\int_Q (\F^{\delta}-\F): \nabla \vp \d \x \d t
$$
for all $\vp \in L^r(0,T; V_r)$. Thus, for $Q_0=I_0 \times B_0$ introduced in Lemma \ref{BDS}, we consider \eqref{dnes} with $\vp \in\C^\infty (Q)$ satisfying $\diver \vp=0$ in $Q_0$ and having the compact support in $Q_0$. Then the boundary term vanishes and we obtain
\begin{align*}
\int_{Q_0}  (\vd -\vv)\cdot \pt \vp \d \x \d t = \int_{Q_0} \left((\Sd -\S) + \vv \otimes \vv - (\vd \otimes \vd) \phi_\delta(|\vd|^2)+ \F^{\delta}-\F\right)\!:\! \nabla \vp \d \x \d t.
\end{align*}
For the sake of ease of exposition, let us denote
\begin{equation}\label{uGG}
\begin{aligned}
\u^\delta &:= \vd - \vv, \\
\H^\delta_1 &:= \Sd - \S, \\
\H^\delta_2 &:= \vv \otimes \vv - (\vd \otimes \vd) \phi_\delta(|\vd|^2)+\F^{\delta}-\F.
\end{aligned}
\end{equation}
Then the triplet $(\u^\delta, \H^\delta_1, \H^\delta_2)$ defined in \eqref{uGG} satisfies the assumptions of Lemma~\ref{BDS}. Recall also that,  due to Lemma~\ref{pidilema}, for $\D\in L^r(Q)$ there is an $\tts\in L^{r'}(Q)$ such that $\G(\tts, \D)=\0$ a.e. in~$Q$ (and, in particular, also in $Q_0$). Referring back to Lemma~\ref{BDS}, we set $\oS:=\S -\tts$ in~\eqref{hvi} and conclude that
\begin{equation}\label{star}
\begin{aligned}
\limsup_{\delta \to 0_+} &\left| \int_{Q_0} (\Sd - \tts):(\Dd - \D) \xi \chi_{Q\setminus Q_{\delta,k}} \right| \\
&= \limsup_{\delta \to 0_+} \left| \int_{Q_0} (\H^\delta_1 + \oS): \nabla \u^\delta \xi \chi_{Q_0\setminus Q_{\delta,k}} \right|\leq C 2^{\frac{-k}{r}}.
\end{aligned}
\end{equation}
Due to \eqref{1816}, $\xi \geq \chi_{\frac{1}{8}Q_0}$; since $\G(\tts,\D)=\0$ and $\G(\Sd, \Dd)=\0$ a.e. in $Q$, the product in the first integral of \eqref{star} is non-negative thanks to~\eqref{G2*} (see Lemma~\ref{pidilema}), and we have
\begin{equation}\label{limhvi}
\limsup_{\delta \to 0_+} \int_{\frac{1}{8}Q_0} \left| (\Sd - \tts): (\Dd - \D) \right| \chi_{Q\setminus Q_{\delta,k}}\d \x \d t \leq C 2^{\frac{-k}{r}}.
\end{equation}
Then for any $a\in(0,1)$,  the following holds:
\begin{align*}
\int_{\frac{1}{8}Q_0} &\left| (\Sd - \tts): (\Dd - \D) \right|^a \d \x \d t \\
&= \int_{\frac{1}{8}Q_0} \left| (\Sd - \tts): (\Dd - \D) \right|^a \chi_{Q_{\delta,k}}\d \x \d t\\
&+ \int_{\frac{1}{8}Q_0} \left| (\Sd - \tts): (\Dd - \D) \right|^a \chi_{Q\setminus Q_{\delta,k}}\d \x \d t \\
&\leq \left(\int_{\frac{1}{8}Q_0} \left| (\Sd - \tts): (\Dd - \D) \right| \chi_{Q_{\delta,k}}\d \x \d t \right)^a |Q_{\delta,k}|^{1-a}\\
&+ \left(\int_{\frac{1}{8}Q_0} \left| (\Sd - \tts):(\Dd - \D) \right| \chi_{Q\setminus Q_{\delta,k}}\d \x \d t \right)^a |Q|^{1-a}\\
&\leq C |Q_{\delta,k}|^{1-a} + C \left(\int_{\frac{1}{8}Q_0} \left| (\Sd - \tts): (\Dd - \D) \right| \chi_{Q\setminus Q_{\delta,k}}\d \x \d t \right)^a.
\end{align*}
Then, as $k \to \infty$, using \eqref{qnk} and \eqref{limhvi}, we conclude
\begin{equation*}
\limsup_{\delta \to 0_+} \int_{\frac{1}{8}Q_0} \left| (\Sd - \tts): (\Dd - \D) \right|^a \d \x \d t \leq C 2^{\frac{-k}{r}} \to 0 \textrm{ as } k\to \infty.
\end{equation*}
However, then also
\begin{equation*}
\left| (\Sd - \tts): (\Dd - \D) \right|^a \to 0 \text{ strongly in } L^1\left(\frac{1}{8}Q_0\right) \textrm{ as } \delta \to 0_+.
\end{equation*}
Due to Egoroff's theorem, for every $\eta$ there exists $Q_{\eta}$ such that $|\frac{1}{8}Q_0 \setminus Q_{\eta}|\leq \eta$, and
\begin{equation*}
\left| (\Sd - \tts): (\Dd - \D) \right|^a \to 0 \text{ strongly in } L^\infty(Q_{\eta}).
\end{equation*}
Consequently,
\begin{equation}\label{CCC1}
(\Sd - \tts): (\Dd - \D) \to 0 \text{ in } L^\infty(Q_{\eta}).
\end{equation}
Since $\lim_{\delta \to 0_+} \int_{Q_{\eta}} \tts: (\Dd - \D)\d \x \d t = 0 $, which follows from \eqref{vvvr}, we obtain  from \eqref{CCC1} that
\begin{equation*}
\lim_{\delta \to 0_+} \int_{Q_{\eta}}\Sd : (\Dd - \D)\d \x \d t = 0,
\end{equation*}
which finally implies (using the weak convergence result for $\Sd$, see \eqref{Sconv}), that
\begin{equation*}
\lim_{\delta \to 0_+} \int_{Q_{\eta}}\Sd : \Dd \d \x \d t = \int_{Q_{\eta}}\S : \D \d \x \d t.
\end{equation*}
According to Lemma~\ref{convGg}, $\G(\S, \D)=\0$ a.e. in $Q_{\eta}$, and we can proceed with $\eta \to 0_+$ to obtain the identification a.e. in $\frac{1}{8} Q_0$. Also, we have that for all $\eta$
\begin{equation}\label{slabasd1}
\Sd : \Dd \tow \S : \D \text{ weakly in } L^1(Q_{\eta}).
\end{equation}

\subsubsection{Energy inequality}
Next, we show that \eqref{enineqr} holds true.
For $0<\er \ll 1$ and $t \in (0,T-\er)$, let $\eta \in \C^{0,1}([0,T])$ be defined as a piece-wise linear function of three parameters, such that
\begin{equation}\label{eta}
\eta(\tau) =
\begin{cases}
1 &\text{if } \tau \in [0,t),   \\
1 + \frac{t-\tau}{\er} &\text{if } \tau \in [t,t+\er),\\
0 &\text{if } \tau \in [t+\er, T].
\end{cases}
\end{equation}
We multiply~\eqref{priprava} by $\eta$, and integrate the result over $(0,T)$ to deduce, after integrating by parts, that
\begin{align*}
\frac{1}{2\er} \int_t^{t+\er} \|\vd(\tau)\|_2^2 \d \tau &+
 \int_{Q^{t+\er}} \!\Sd:\Dd \eta \d \x \d \tau +  \int_{\Gamma^{t+\er}} \! \fidvd \sd \cdot \vd \eta \d S \d \tau \\
&= \int_0^{t+\er} \langle \b^{\delta}, \vd \eta \rangle_{V_r} \d \tau + \frac12 \|\vv_0\|_2^2.
\end{align*}
The next step is to take the limit as $\delta \to 0_+$. We know that  $(\Sd:\Dd)\geq 0$ and $(\sd \cdot \vd)\geq 0$ because $\G(\Sd,\Dd)=\0$ in $Q$, $\g(\sd, \vd)=\0$ on $\Gamma$ and we have~\eqref{G2*} from Lemma~\ref{pidilema} and an analogous result holds for the function $\g$ as well. Therefore, from the above identity we deduce
\begin{equation}\label{iuge}
    \begin{aligned}
\frac{1}{2\er} \int_t^{t+\er} \|\vd(\tau)\|_2^2 \d \tau &+ \int_{Q^t \cap Q_{\eta}} \Sd:\Dd \d \x \d \tau + \int_{\Gamma^t \cap \Gamma_{\eta}} \fidvd \sd \cdot \vd \d S \d \tau \\
&\leq \int_0^{t+\er} \langle \b^{\delta}, \vd \eta \rangle_{V_r} \d \tau + \frac12 \|\vv_0\|_2^2.
\end{aligned}
\end{equation}
For the first term, we can use the weak lower semicontinuity of the norm. For the products $\Sd:\Dd$ and $\fidvd\sd \cdot \vd$, we use~\eqref{slabasd1},~\eqref{slabasv1} and the weak convergence of $\vd$ on $\Gamma_{\eta}$.
For the duality term, we use \eqref{bstro} and  \eqref{vvvr}, and get that
\begin{align*}
\frac{1}{2\er} \int_t^{t+\er} \|\vv(\tau)\|_2^2 \d \tau &+ \int_{Q^t\cap Q_{\eta}} \S:\D \d \x \d \tau + \int_{\Gamma^t \cap \Gamma_{\eta}} \s\cdot\vv \d S \d \tau \\
&\leq \int_0^{t+\er} \langle \b, \vv \eta \rangle_{V_r} \d \tau + \frac12 \|\vv_0\|_2^2.
\end{align*}
Next, we proceed with $\er, \eta \to 0_+$, then $Q^t\cap Q_{\eta} \to Q^t$ and $\Gamma^t \cap \Gamma_{\eta} \to \Gamma^t$, and finally, thanks to $\vv \in \C_w([0,T];H)$ and the fact that the other terms are well-defined, we obtain the energy inequality \eqref{enineqr} for any $t \in (0,T)$.

\subsection{Attainment of the initial condition}

Considering $\eta$ introduced in~\eqref{eta}, we multiply \eqref{WFd} by $\vp \eta$, where $\vp \in V_z$ is arbitrary. Integrating the result  over $(0,T)$, we get
\begin{align*}
\iT \langle\pt \vd, \vp\rangle_{V_r} \eta \d \tau &- \iq \left((\vd \otimes \vd) \phi_\delta(|\vd|^2)\right) : \nabla \vp \eta \d \x  \d \tau \\
&+ \iq \Sd : \DD \vp \eta\d \x \d \tau + \ig \fidvd \sd \cdot \vp \eta \d S \d \tau = \iT \langle \b^{\delta}, \vp \eta \rangle_{V_r} \d \tau.
\end{align*}
Next, we apply integration by parts in the first term (using the properties of $\eta$ and the fact that $\vp$ is independent of $t$).  Then we take the limit $\delta \to 0_+$. Using arguments from Sect.~\ref{subposl}, in paricular the convergence result~\eqref{vdovd} to take the limit in the convective term, we conclude that
\begin{align*}
\frac{1}{\er} \int_t^{t+\er} \io \vv \cdot  \vp \d \x \d \tau &- \int_{Q^{t+\er}} (\vv \otimes \vv) : \nabla \vp \eta \d \x  \d \tau + \int_{Q^{t+\er}} \S:\DD \vp \eta \d \x \d \tau \\
&+ \int_{\Gamma^{t+\er}} \s\cdot\vp \eta \d S \d \tau = \int_0^{t+\er} \langle \b, \vp \eta \rangle_{V_r} \d \tau + \io \vv_0 \cdot \vp \d \x.
\end{align*}
Since $\vv \in \C_w([0,T];H)$ and $\vp$ is independent of time, we can let $\er \to 0_+$ to conclude
\begin{align*}
\io \vv(t) \cdot  \vp \d \x  &- \int_{Q^t} (\vv \otimes \vv) : \nabla \vp \d \x  \d \tau + \int_{Q^t} \S:\DD \vp \d \x \d \tau \\
&+ \int_{\Gamma^t} \s\cdot\vp \d S \d \tau = \int_0^{t} \langle \b, \vp  \rangle_{V_r} \d \tau + \io \vv_0 \cdot \vp \d \x.
\end{align*}
Now, standard density arguments imply that
\begin{equation*}
\vv(t) \tow \vv_0 ~\text{ weakly in }~L^2(\o).
\end{equation*}
Also, taking the limes superior for $t\to 0_+$ in the energy inequality \eqref{enineqr}, we obtain that $\limsup_{t \to 0_+} \|\vv(t)\|_2^2 \leq \|\vv_0\|_2^2$. The last two pieces of information imply the strong convergence in $H$ as claimed in \eqref{IC}.

\subsubsection{Existence of an integrable pressure}

In order to reconstruct the pressure, we need to assume that $\o \in \C^{1,1}$. The procedure to obtain an integrable pressure for problems with slipping boundary conditions, is explained in~\cite{BGMS2} in detail. This is why we merely show here a formal estimate concerning the ``best" integrability of $p$. To do so, we assume that there exists an integrable pressure~$p$ such that~\eqref{WFp} holds. Moreover, we may assume that~$\io p(t) \d \x =0$ for a.a.~$t$. Next, we find~$\phi$ solving the equation
\begin{align*}
\Delta \phi &= |p|^{z'-2} p - \frac{1}{|\o|} \io |p|^{z'-2} p \d \x &&\text{in } \o, \\
\nabla \phi \cdot \n &= 0 &&\text{on } \po.
\end{align*}
Using classical theory we know that such a $\phi$ exists and satisfies the estimate
\begin{equation}\label{phiz}
    \|\phi\|_{W^{2,z}(\o)}^z \leq C \io |p|^{z'} \d \x.
\end{equation}
Setting $\vp:= \nabla \phi$ in~\eqref{WFp}, we have (note that $\vp \in \cV_z$ and also that the term with the time derivative vanishes)
\begin{equation}
 \io |p|^{z'} \d \x = - \int_\o (\vov):\nabla^2 \phi \d \x + \int_\o \S:(\DD \nabla \phi) \d \x + \int_{\po} \s \cdot \nabla \phi \d S - \langle \b, \nabla \phi \rangle_{\cV_z}.
\end{equation}
By use of~\eqref{phiz}, the H\"older inequality and the Trace Theorem we can deduce that
\begin{equation*}
    \|p\|_{z'}^{z'} \leq C \left( \|\vov\|_{z'} + \|\S\|_{z'}  + \|\s\|_{L^{z'}(\po)} + \|\b\|_{\cV_z^*} \right)^{z'}.
\end{equation*}
After integrating it over $(0,T)$, thanks to the definition of $z'$ and a stronger assumption on $\b \in L^{z'}(0,T;\cV_z^*)$, all terms on the right-hand side are bounded and we conclude that $p\in L^{z'}(0,T;L^{z'}(\o))$.

\section{Extensions: existence results for related problems (a summary)}\label{Sect7}

The above established result concerns isothermal homogeneous and incompressible fluids; these properties can be seen as limitations and one may wishes to develop a theory in the similar spirit as above for heat-conducting or inhomogeneous or compressible fluids or for fluids that share more of these properties. Also the constitutive equation \eqref{pepa9} does not cover viscoelastic rate-type or integral models. Similarly, the boundary condition \eqref{pepa10} does not include dynamic boundary conditions. Below, we provide references that can be relevant to anyone who would like to extend the study in the directions indicated.

When one wishes to include the dependence of the viscosity on the temperature, the system of governing equations has to be completed by the formulations of the balance of energy and the second law of thermodynamics. One also needs to specify the constitutive equation relating the heat flux to the temperature gradient, boundary conditions for the temperature etc. These extensions give rise, within the context of weak solutions, to several concepts of solution. %Here we mention our system can be completed also by the equation for the energy to have the generalized incompressible Navier--Stokes--Fourier system.
A sound long-time and large-data existence theory for \emph{heat-conducting} fluids described by incompressible Navier--Stokes--Fourier equations goes back to \cite{FeMa06,BuMaFe09}. The first existence result for nonlinear models of the power-law type is due to Consiglieri~\cite{Luisa} for $r\geq \frac{3d+2}{d+2}$. Then, a very similar theory was obtained for smaller $r$'s in \cite{BuMaRa09} but merely for explicit models of the type $\S= 2\nu(\theta, |\D|^2)\D$, where $\theta$ is the temperature. The first fully implicit approach with material parameters depending on the temperature was  developed in \cite{MaZa18}, where the authors dealt with a specific activated model with the activation depending on the temperature that also allowed the model to range from activated Euler through the Navier--Stokes regime to a Bingham-type response. %Based on the results mentioned above, it seems to be clear, that our fully implicit setting described by the conditions (G1)--(G4) can be easily adapted into the heat conducting case.}

The mathematical theory for unsteady flows of \emph{inhomogeneous} Navier-Stokes fluids is developed in the book by P.-L. Lions~\cite{pll96}. Inhomogeneous isothermal fluids of power-law type are analyzed in \cite{FrRu10}, while heat-conducting processes for such fluids are treated in \cite{FrMaRu}. \emph{Compressible} non-Newtonian fluids of power-law type serve represent a completely open field from the point of view of the large-data analysis of relevant initial- and boundary-value problems (see \cite{AbFeNo} for one of the first attempts and further references).

Regarding \emph{viscoelastic rate-type} fluids, where for a part of the Cauchy stress one has an additional evolutionary (nonlinear) equation, we distinguish two basic classes: without stress diffusion and with stress diffusion. The first class includes the standard Maxwell, Oldroyd-B or Giesekus models, and the additional equation is of transport type. The long-time and large-data mathematical theory goes back to \cite{LiMa00}, where Oldroyd-B type models with a corotational time derivative are studied. The type of derivative is however non-physical and simplifies the analysis tremendously. For Giesekus  type of models with more general objective derivatives the idea of existence proof is presented in \cite{Ma11} and rigorously proven in the planar case in \cite{BuLoLuMa22}. Regarding the analysis of viescoelastic rate-type fluids with stress diffusion, we refer to \cite{BaBuMa21} where a very robust theory is developed and other relevant studies are cited. %. Similarly as for Navier--Stokes--Fourier system, we strongly believe that one can directly adapt the method presented in this paper to cover also more general models of viscoelastic fluids mentioned above.

Finally, \emph{dynamic} boundary conditions  (see
\cite{Hatzi} for their relevance to observations connected with the experiments regarding molten polymers) are the subject of a recent investigation from the point of view of analysis of partial differential equations. The theory for the Stokes system with such dynamic boundary conditions is developed in \cite{AbBuMa21}.

Extensions in other directions are possible. Examples include the analysis of rapidly shear-thickening fluids (see \cite{GwSw08}) or fluids with a~priori bounded velocity gradient (see \cite{MiRo2016, BuHrMa}). Also, one can consider instead of \eqref{pepa9} a more general class of incompressible fluids given by the relation
$\G(\T,\D) = \0$, which allows one to naturally include naturally fluids with pressure and shear-rate dependent viscosity; here we refer to \cite{BuMaRa07,BuMaRa09,BuMaRa09b,BuZa15} for further details.

%\begin{itemize}
%\item temperature - Navier--Stokes--Fourier %\cite{FeMa06,BuMaFe09,MaZa18} \cite{AF} %\cite{BuMaRa09}
%\item dynamic boundary conditions \cite{AbBuMa21}
%\item viscoelasticity - with stress diffusion / %without stress diffusion \cite{BaBuMa21}, %\cite{BuLoLuMa22} \cite{LiMa00} \cite{Ma11} %\cite{KrPoSa15}
%\end{itemize}

\section{Uniqueness, smoothness, open problems and concluding remarks}\label{Sect8}

 We have studied long-time and large-data mathematical properties of unsteady internal flows of incompressible fluids with frictional properties characterized by implicit constitutive equations in the bulk and on the boundary. The developed theory that has origin in the seminal works of O.~A.~Ladyzhenskaya addresses positively the question concerning existence of weak solutions to large classes of fluids as well as boundary conditions. The structural assumptions characterizing the admissible class of fluids and boundary conditions are expressed in terms of basic tools of calculus and can be checked directly for a given constitutive equation without any deeper knowledge of concepts of the operator theory. Despite broad applicability of the developed theory to many models in different scientific areas, the analysis for fluids satisfying~(G1) and~(G4) but having  \emph{non-monotone} response (as for example the model computationally tested in \cite{Janecka2019}) is an open problem. (Note however that non-monotone responses in the boundary conditions can be included, see \cite{BuMa19} for details.)

The main achievement of this study lies in a novel existence theory. Ladyzhenskaya's interest in these fluids was however motivated by the uniqueness and smoothness of these solutions. Concerning the uniqueness of weak solution, one can easily observe that for models fulfilling (G1)--(G4) we automatically/easily obtain the uniqueness of the velocity field for all $r>1$ provided that we neglect the convective term $\diver(\vv\otimes\vv)$; compare with \cite{BMM21}. For the complete model, i.e. for model with the convective term, already Ladyzhenskaya was able to show the uniqueness of a weak solution for the models of the form
\begin{equation}
  \S= 2\nu_*\D + 2\tilde\nu_* |\D|^{r-2}\D \quad \textrm{ or } \quad \S=2\nu_* (1+|\D|^2)^{\frac{r-2}{2}}\D \label{pepa20}
\end{equation}
provided that $r\geq \frac{d+2}{2}$. To date the most general uniqueness result is due to \cite{BuKaPr19}, where the authors prove uniqueness of a weak solution in three dimensions (for sufficiently regular data) and a class of models having the $r$-growth with $r\geq \frac{11}{5}$ and satisfying, for a certain $C_* >0$
\begin{equation*}
\begin{split}
    (\S_1 - \S_2):&(\DD_1 - \DD_2) \geq C_* \left(|\DD_1 - \DD_2|^2 +|\DD_1 - \DD_2|^r\right)\\
    &\text{ for all } (\S_i, \DD_i) \text{ such that } \G(\S_i,\DD_i)=\0,~i=1,2,
\end{split}
\end{equation*}
which is fulfilled by the models given in \eqref{pepa20} above, but it is much more restrictive in comparison with the condition (G2*) (or its equivalent form (G2)) needed in the existence theory.

We can slightly strengthen these uniqueness results by considering fluids that behave as the Navier-Stokes fluid prior the activation, i.e. for  $|\D|\leq \delta_*$, where $\delta_*>0$ can be arbitrary, and behave as a power-law fluid with $r\geq \frac{3d+2}{d+2}$ once the activation takes place, i.e. $|\D|> \delta_*$. Mathematically, such a model is described by the constitutive equation of the form
$$
\S=2\nu_* \D + 2\tilde\nu_*\frac{(|\D|-\delta_*)_+}{|\D|}|\D|^{r-2}\D.
$$
For this model, one can then establish both the existence and the uniqueness of a weak solution for $r$ sufficiently large ($r\geq\frac{3d+2}{d+2}$). We wish to emphasize that the constant $\delta_*$ can be chosen arbitrarily large. On the other hand, for $r<\frac{3d+2}{d+2}$ we have a counterexample to uniqueness thanks to \cite{MR4328053} in the class of very weak solutions. Hence, a natural open problem is the (non)uniqueness of a weak solution for smooth data in natural function spaces also for $r<\frac{11}{5}$ in dimension three.

Finally, concerning smoothness of weak solutions to nonlinear models studied in this paper, there is one striking open problem. Independently of whether one excludes or includes the convective term and independently of the value of the parameter $r$, it is not clear (even in the situation when we know that there is a unique weak solution) whether, for smooth but large data, there exists a global-in-time $C^{1,\alpha}$-solution for any special case of the problem \eqref{pepa3}, \eqref{pepa9} and \eqref{pepa10} in three dimensional setting. The regularity theory in two dimensions is available, see e.g. \cite{BuKaPr19,DiKaSc14,KaMASt98,KaMaSt02,BuMaSh14}, the theory in dimension three is however basically untouched.

%\bibliographystyle{amsplain}
%\bibliography{biblio}

\providecommand{\bysame}{\leavevmode\hbox to3em{\hrulefill}\thinspace}
\providecommand{\MR}{\relax\ifhmode\unskip\space\fi MR }
\providecommand{\MRhref}[2]{%
  \href{http://www.ams.org/mathscinet-getitem?mr=#1}{#2}
}
\providecommand{\href}[2]{#2}

\end{document}